\documentclass[journal]{IEEEtran}
\usepackage{cite} 

\usepackage{graphicx} 
\usepackage[dvipsnames]{xcolor}  
\usepackage{mdframed}
\usepackage{subfigure}

\usepackage{amsmath,amssymb,mathrsfs}
\usepackage{mathtools,xparse}
\usepackage[shortlabels]{enumitem} 
\newcommand{\CCC}{{\mathbf C}}
\newcommand{\PPP}{{\mathbf P}}
\DeclarePairedDelimiter{\abs}{\lvert}{\rvert}
\DeclarePairedDelimiter{\norm}{\lVert}{\rVert}
\newcommand{\rfb}[1]{\mbox{\rm
   (\ref{#1})}\ifx\undefined\stillediting\else:\fbox{$#1$}\fi} 
\newcommand{\rarrow} {\mathop{\rightarrow}}                
\newcommand{\FORALL} {{\hbox{$\hskip 11mm \forall \;$}}}
\newcommand{\e}      {{\varepsilon}}
\newcommand{\m}      {{\hbox{\hskip 1pt}}}
\newcommand{\mm}     {{\hbox{\hskip 0.5pt}}}
\newcommand{\nm}     {{\hbox{\hskip -3pt}}}
\renewcommand{\l}    {{\lambda}}
\renewcommand{\Re}   {{\rm Re\,}}

\newcommand{\rline}  {{\mathbb R}}
\newcommand{\Kscr}   {{\mathscr{K}}}
\newcommand{\Lscr}   {{\mathscr{L}}}

\newcommand{\Sscr}   {{\mathscr{S}}}
\newcommand{\dd}     {{\rm d}\hbox{\hskip 0.5pt}}
\newcommand{\sbluff} {{\hbox{\raise 9pt \hbox{\mm}}}}
\newcommand{\half}   {{\frac{1}{2}}}
\newcommand{\bbm}[1]{\left[\begin{matrix} #1 \end{matrix}\right]}
\newcommand{\sbm}[1]{\left[\begin{smallmatrix} #1
	\end{smallmatrix}\right]}

\newtheorem{theorem}{Theorem}[section]

\newtheorem{lemma}[theorem]{Lemma}
\newtheorem{proposition}[theorem]{Proposition}

\newtheorem{remark}[theorem]{Remark}
\newtheorem{assumption}{Assumption}

\begin{document}
\title{Saturating PI Control of Stable Nonlinear Systems Using 
       Singular Perturbations}

\author{Pietro~Lorenzetti and~George~Weiss
\thanks{The authors are team members in the ITN network ConFlex, 
   funded by the European Union's Horizon 2020 research and innovation
   program under the Marie Sklodowska-Curie grant agreement no. 
   765579.}%
\thanks{P. Lorenzetti is with the Department of Electrical 
   Engineering Systems, Faculty of Engineering, Tel Aviv University,
   Ramat Aviv 69978, Israel (e-mail: plorenzetti@tauex.tau.ac.il)}%
\thanks{G. Weiss is with the Department of Electrical Engineering
   Systems, Faculty of Engineering, Tel Aviv University, Ramat Aviv 
   69978, Israel (e-mail: gweiss@tauex.tau.ac.il)}}

\maketitle

\begin{abstract} 
This paper presents an anti-windup PI controller, using a
saturating integrator, for a single-input single-output stable
nonlinear plant, whose steady-state input-output map is increasing. We
prove that, under reasonable assumptions, there exists an upper bound
on the controller gain such that for any constant reference input, the
corresponding equilibrium point of the closed-loop system is
exponentially stable, with a ``large'' region of attraction. When the
state of the closed-loop system converges to this equilibrium point,
then the tracking error tends to zero. The closed-loop stability
analysis employs Lyapunov methods in the framework of singular
perturbations theory. Finally, we show that if the plant satisfies the
asymptotic gain property, then the closed-loop system is globally
asymptotically stable for any sufficiently small controller gain. The
effectiveness of the proposed PI controller is proved by showing how
it performs as part of the control algorithm of a synchronverter (a
special type of DC to AC power converter).
\end{abstract}

\IEEEpeerreviewmaketitle

\section{Introduction} 

\IEEEPARstart{P}{Proportional-integral} (PI) control is extensively
used to achieve robust asymptotic tracking and disturbance rejection
of constant external signals (it is the simplest instance of the
internal model principle). When the plant is an uncertain linear
system, closed-loop stability can be achieved if the control gains are
sufficiently small and the plant fulfils certain stability conditions
\cite{Fliegner2003,Morari1985}. The theory has been extended to
nonlinear systems, where global results have been derived in
\cite{Desoer1985} for PI control and in \cite{Simpson-Porco2020}
for integral control, and local ones in \cite{Isidori1990}
for the more general internal model principle.  Later, regional and
semiglobal results have been presented in \cite{Isidori1997},
\cite{Khalil2000} for specific classes of nonlinear systems using 
high-gain observers.
For extensions of PI control to infinite dimensional linear systems
see \cite{Hamalainen1996,Logemann1999,Rebarber2003}.

In the presence of actuator limitations, e.g. saturation, there can be
a mismatch between the actual plant input and the controller output.
When this happens, the controller output does not drive the plant and,
as a result, the states of the controller are wrongly updated. If the
controller is unstable (such as PI), then the state of the controller
may reach a region far from the normal operating range, an effect
called {\em controller windup} \cite{Kothare1994}. (Windup can also be
caused by a temporary fault in the system.) Windup may cause long
transients, oscillations and even instability. To prevent windup,
several anti-windup techniques have been investigated in the last 50
years. Mainly, two design approaches can be identified (see
\cite{Tarbouriech2009}): {\em one-step} design, where both the
actuator limitations and the performances are taken into account
during the controller design, and {\em two-steps} design (also known
as {\em anti-windup compensation}), where first a nominal controller
is designed ignoring the actuator nonlinearities, focusing on the
performance, and then an {\em anti-windup compensator} is added to
handle the actuator nonlinearities so that the closed-loop performance
degrades ``gracefully'' when actuator limitations occur.

In the field of linear systems, one of the first systematic
anti-windup methods is the {\em conditioning technique} introduced in
\cite{Hanus1987} as an extension of the \textit{back calculation
method} presented in \cite{Fertik1967}. Later, this and other designs
have been collected in \cite{Astrom1989}, where the more general
\textit{observer approach} is formulated for a class of anti-windup
compensation schemes, which includes the conditioning technique. For
the application of internal model control to the anti-windup problem,
see \cite{Zheng1994}. Even though this method does not excel in terms
of performance, it has been proved later (see \cite{Turner2007}) to be
the optimal one in terms of robustness.  One of the first theoretical
frameworks for the study of two-steps anti-windup design for LTI
systems is in \cite{Kothare1994}, where a general coprime-factor
scheme is introduced, unifying all the known anti-windup
methods. Related later work is in \cite{Edwards1998}, where a generic
approach in the form of an $H_\infty$ optimization framework is
proposed. A comprehensive review of anti-windup techniques specific
for PID controllers is \cite{Peng1996}. Additional results on
anti-windup schemes for PI and PID controllers are respectively in
\cite{Choi2009}, where the conditional integrator scheme is improved
by setting a specific initial value for the integrator when switching
from P to PI control mode, and in \cite{Hodel2001}, where a variable
structure PID is presented, able to switch between normal operating
mode and anti-windup compensator mode.  A recent result on anti-windup
low-gain integral control for multivariable linear systems subject to
input nonlinearities is in \cite{Guiver2017}. A treatment of the
$\Lscr_2$ anti-windup problem is in \cite{Zaccarian2002}. In recent
years, more attention has been devoted to LMIs to tackle complex
anti-windup scenarios. In fact, the actuator nonlinearities can be
modelled as sector nonlinearities so that the anti-windup controller
design is recast into a convex optimization problem under LMI
constraints \cite{Tarbouriech2009}. In this context, in
\cite{Turner2007} the uncertainty in the system model is taken into
account, and in \cite{Wu2004} regional stability is guaranteed in the
presence of input saturation for an exponentially unstable linear
system.
These and other LMI-based methods can be found in the comprehensive
surveys \cite{Hippe2006,Tarbouriech2011,Tarbouriech2009}.

Nonlinear systems have received less attention in the anti-windup
literature. The first notable results in this direction are based on
input-output linearization techniques, combined with internal model
control, for nonlinear affine systems. Among them, we refer to 
\cite{Hu2000} for single-input-single-output nonlinear systems and to
\cite{Kapoor1999,Kendi1997} for multivariable systems.
For Euler-Lagrange systems, an anti-windup scheme inducing global
asymptotic stability and local exponential stability is proposed in
\cite{Morabito2004}. For LMI-based optimization techniques, the reader
is referred to \cite{daSilva2012}, \cite{Rehan2013}, and the
references therein.

A one-step anti-windup design for nonlinear systems has been recently
proposed in \cite{Konstantopoulos2016}, where a \textit{bounded
integral controller} (BIC) is presented. Under suitable assumptions
(input to state practical stability), the BIC generates a bounded
control signal, guaranteeing boundedness of the plant state
trajectory, while achieving tracking for constant references.
According to \cite{Konstantopoulos2016}, most integral control
algorithms for nonlinear plants require detailed knowledge of the
plant, leading to high complexity.  Therefore, a simple anti-windup
controller that allows a rigorous stability analysis, based on simple
assumptions about the plant, is of significance.

This paper proposes a one-step anti-windup PI design for stable
nonlinear systems. This PI controller uses a \textit{saturating
integrator}, which ensures that the integrator state is constrained to
a compact interval chosen based on physical constraints. The
saturating integrator is straightforward to implement and it has
proved to be effective in power electronics applications, see for
instance \cite{Kustanovich2020b}, \cite{Lifshitz2015} and
\cite{Natarajan2017}. Our main result, Theorem \ref{thm:stab}, shows
that under reasonable assumptions, for small enough controller gain,
the closed-loop system is (locally) exponentially stable and the
region of attraction of its equilibrium point contains a curve of
equilibrium points of the open-loop system.

Note that results closely related to our Theorem \ref{thm:stab},
considering only integral control, with many proofs missing, and
without using singular perturbation theory, have been presented in the
conference paper \cite{Weiss2016}. We believe that the new singular
perturbations approach is much more elegant, and more suitable for
generalizations. A preliminary and partial presentation of our
results, without applications, is in our recent conference paper
\cite{Lorenzetti2020}, again considering only integral control.
Moreover, the global stability results from Section \ref{sec:5} here
are substantially stronger than those in \cite{Lorenzetti2020}.

The paper is organized as follows. In Section \ref{sec2} the
saturating integrator is introduced and the control system under
consideration is described in precise terms. In Section \ref{sec3} the
closed-loop system equations are rewritten as a singular perturbation
model. In Section \ref{sec:stability_analysis} the main result is
presented, with the stability analysis of the closed-loop system using
singular perturbation methods. In Section \ref{sec:5} we assume that
the plant satisfies the asymptotic gain property around each
equilibrium point. Then, we show that for small enough controller
gain, the closed-loop system is globally asymptotically stable.
Finally, in Section \ref{sec:ex} we illustrate our theory using two
examples: a toy example and a real application in the control of a
certain type of three phase inverters.

\section{The closed-loop system} \label{sec2} 

We describe the control system that will be investigated here. The 
nonlinear plant $\PPP$ to be controlled is described by:
\begin{equation} \label{eq:P}
   \dot{x} \m=\m f(x,u), \qquad y \m=\m g(x),
\end{equation}
where $f\in C^2(\rline^n\times\rline;\rline^n)$ and $g\in {\rm Lip}
(\rline^n;\rline)$. In this paper, Lip denotes a space of locally
Lipschitz functions, that is, functions that are Lipschitz on any
compact set. 

The control objective is to make the output signal $y$ track a
constant reference signal $r\in Y\vcentcolon=[y_{min},y_{max}]\subset
\rline$, using an input signal that in steady state takes values in
the range $U\vcentcolon=[u_{min},u_{max}] \subset\mathbb{R}$ (here
$u_{min}<u_{max}$). In order to achieve this goal, a type of
anti-windup integral controller is used, which we call the {\em
saturating integrator}. We define the positive (negative) part of a
real number $w$ by $w^+=\max\{w,0\}$ ($w^-=\min\{w,0\}$), so that
$w=w^++w^-$. The {\em saturating integrator} is a system with input
$w$ and state $u_I$, described by \vspace{-3mm}
\begin{equation} \label{eq:sat}
   \begin{gathered} \dot{u}_I \m=\m \Sscr(u_I,w), \\
   \m\hspace{-3mm}\text{where} \qquad \Sscr(u_I,w) \m=\m \begin{cases} 
   w^+ & \text{if} \quad u_I\leq u_{min}, \\ w & \text{if} \quad u_I 
   \in (u_{min},u_{max}), \\ w^- & \text{if} \quad u_I\geq u_{max}. 
   \end{cases} \end{gathered}
\end{equation}
The idea behind \rfb{eq:sat} is that if $u_I\geq u_{max}$ ($u_I\leq
u_{min}$), we do not allow $u_I$ to move further into the forbidden
region.  In this way, assuming that $u_I(0)\in U$, the state
trajectory $u_I(t)$ is constrained in $U$ for all $t\geq0$.
A similar controller is described in \cite{Hodel2001}, under
the name of \textit{conditionally freeze integrator}, however no
closed-loop stability results are given. We also mention that
the \textit{limited integrator} present in the MATLAB/Simulink software 
is of the form \rfb{eq:sat}.

\begin{remark}
As discussed in \cite{Biannic2009}, the saturating integrator \rfb{eq:sat} can
be approximated by the nonlinear system with input $w$, state $v_I$ and output
$\bar{v}_I$ described by \vspace{-1mm} $$\dot{v}_I=w-\lambda(v_I-\bar{v}_I),
\qquad \bar{v}_I={\rm sat}\m v_I,\vspace{-1mm}$$ where ${\rm sat}$ is the
saturating function with lower bound $u_{min}$ and upper bound $u_{max}$, and
the gain $\lambda>0$ is sufficiently large. More precisely, as stated in
\cite[Lemma~5.1]{Biannic2009}, the approximation error
$e_I\vcentcolon=\bar{v}_I-u_I$, where $u_I$ is the state of \rfb{eq:sat}, tends
to zero for increasing values of $\lambda$. The advantage of this representation
is that, in case of a linear plant $\PPP$, the stability analysis can be
performed using the LMI-based tools developed, for instance, in
\cite{Biannic2009,Sofrony2009,Tarbouriech2009,Tarbouriech2011}.
\end{remark}

\begin{figure} \centering 
   \includegraphics[width=0.38\textwidth]{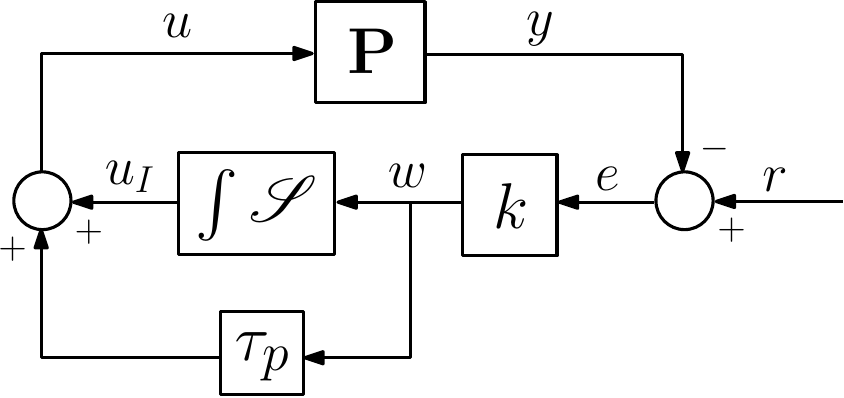}
   \caption{Closed-loop system formed by the plant $\PPP$, the
   saturating integrator $\int\Sscr$ and the constant gains $k>0$,
   $\tau_p\geq0$, with the reference $r$.} \label{fig:clos}
   \vspace{-3mm} \end{figure}

As shown in Fig. \ref{fig:clos}, we connect the saturating integrator
in parallel with a proportional gain $\tau_p\geq 0$, and the input to
these two blocks comes through a gain $k>0$. If the integrator
does not saturate, then this is a linear PI controller with transfer
function
\begin{equation} \label{eq:C}
   \CCC(s) \m=\m k\left( \tau_p+\frac{1}{s} \right) \m.
\end{equation}
(Here $\tau_p\geq0$ is considered to be given, while
$k>0$ to be chosen according to the discussion in Remark \ref{rmk:bound_kappa}.)
We assume that $r\in Y$. Most of the time we have that $u(t)\in U$,
thanks to the saturating integrator, and thanks to the fact that in
steady state, $e(t)=0$.  During transients the signal $u$ may exit $U$
because of the term $\tau_p w$, but $u_I$ will remain in $U$. The
closed-loop system from Fig. \ref{fig:clos} is given by
\begin{equation} \label{Trump}
   \begin{cases} \begin{gathered}
   \m \ \dot{x} \m=\m f(x,u_I+\tau_p k(r-g(x))), \\
   \dot{u}_I \m=\m \Sscr(u_I,k(r-g(x))), \end{gathered} \end{cases}
\end{equation}
with state $(x,u_I)$ and state space $X\vcentcolon=\rline^n\times U$.

An informal statement of our main result (see Theorem \ref{thm:stab})
is that, under reasonable assumptions on the plant $\PPP$, the
following holds: For any small enough $k>0$, the closed-loop system
\rfb{Trump} shown in Fig. \ref{fig:clos}, with any constant reference
$r\in Y$, is locally exponentially stable around an equilibrium point,
with a ``large'' region of attraction. When the state converges to
this equilibrium point, then the tracking error $e\vcentcolon= r-y$
tends to zero at an exponential rate. Moreover, this result remains
true also if $r$ has finitely many step changes, if the time interval
between consecutive changes is large enough.

A related result, but without limitations on $u$, assuming global
exponential stability of the plant's equilibrium point for any
constant input $u\in \mathbb{R}^m$, has been given in
\cite{Desoer1985}. Recently, the result from \cite{Desoer1985} has 
been generalized in \cite{Simpson-Porco2020}, where the assumption on the monotonicity
of the plant steady-state input-output map has been
replaced by the (weaker) uniform infinitesimal contracting property of the reduced dynamics.

Before addressing the stability analysis of the closed-loop system
\rfb{Trump}, some care has to be taken to define its trajectories.
First, it has to be ensured that the state trajectories $u_I$ of the
saturating integrator \rfb{eq:sat} are well-defined for any input
$w\in L^1[0,t]$. For this we resort to a density argument: Note that
for a polynomial input $w$, the state trajectory $u_I$ is easy to
define. (The key is that a polynomial has finitely many
zeros. Problems with defining the solution of \rfb{eq:sat} may arise
when the zeros of $w$ have an accumulation point, for instance if
$w(t)=t\sin(\frac{1}{t})$ and $u_I(0)=u_{max}$.)

Let $u_1,u_2$ be the state trajectories of \rfb{eq:sat} 
corresponding to the polynomials $w_1,w_2$. Then by an elementary argument
\begin{equation}\label{eq:blair}
	\abs{u_2(t)-u_1(t)}\leq\abs{u_2(0)-u_1(0)}+\int_{0}^{t}
	\abs{w_2(\sigma)-w_1(\sigma)}d\sigma.
\end{equation}
This shows that $u_I(t)$ depends Lipschitz continuously both on
$u_I(0)$ and also on $w$ considered with the $L^1$ norm. For
$u_2(0)=u_1(0)$ we can write the last estimate as
\begin{equation} \label{eq:lipschitz_u_w}
   |u_2(t)-u_1(t)| \leq \|w_2-w_1\|_{L^1[0,t]}.
\end{equation}
Hence, by continuous extension, we can define $u_I(t)$ for any initial
state $u_I(0)$ and for any input $w\in L^1[0,t]$ (because the
polynomials are dense in $L^1[0,t]$). Next, we show the local
existence and uniqueness of the state trajectories of \rfb{Trump},
which is not trivial due to the discontinuity of $\mathscr{S}$.

\textbf{Notation.} For any $\delta>0$, denote $U_\delta\vcentcolon=
[u_{min}-\delta,u_{max}+\delta]$.\vspace{-4mm}
{\color{blue}
\begin{proposition} \label{prop:uniqueness}
Let $\PPP$ be described by \rfb{eq:P} and let $\int\Sscr$ be the
saturating integrator from \rfb{eq:sat}. For every $k,\tau_p\in\mathbb{R}$, 
$x_0\in\rline^n$, $\delta>0$, $u_0\in U_\delta$, and every constant
$r\in\rline$ there exists a $\tau\in(0,\infty]$ such that the
closed-loop system \rfb{Trump} (shown in Fig.  {\rm\ref{fig:clos}})
has a unique state trajectory $(x,u_I)$ defined on $[0,\tau)$, such
that $x(0)=x_0$ and $u_I(0)=u_0$. If $\tau$ is finite and maximal
(i.e., the state trajectory cannot be continued beyond $\tau$), then
$\limsup_{t\rarrow\tau}\|x(t)\|=\infty$.
\end{proposition}}

\vspace{1mm}
For the proof see Appendix \ref{app:A}.

\section{Formulation as a standard singular perturbation model} 
         \label{sec3} 

To perform the stability analysis of the closed-loop system
\rfb{Trump}, the idea is to regard the constant gain $k>0$ as a
``sufficiently small'' parameter such that \rfb{Trump} can be
rewritten as a standard singular perturbation model (see
\cite[Chapt.~1]{Kokotovic1999} or \cite[Chapt.~11]{Khalil2002}). Then,
we can separately analyze the stability of the reduced (slow) model
and of the boundary-layer (fast) system, and, using the tools from
singular perturbation theory, obtain stability results for
\rfb{Trump}. Note that in the less general formulation in
\cite[Sect.~11.5]{Khalil2002}, the functions determining the
closed-loop system are required to be locally Lipschitz, while in
\cite[Chapt.~7]{Kokotovic1999} it is only required that a unique
(local) solution exists for the closed-loop system, which we have
proved in Proposition \ref{prop:uniqueness}.

The following assumption is common in the singular perturbation theory (see
\cite{Desoer1985},\cite{Simpson-Porco2020},
\cite[Chapt.~11]{Khalil2002}, \cite{Kokotovic1999}) and in the theory for
nonlinear systems with slowly varying inputs (see \cite{Kelemen1986},
\cite{Lawrence1990}, \cite[Sect.~9.6]{Khalil2002}).

\begin{mdframed}
\begin{assumption} \label{ass_1}
There exists $\delta>0$ and a function $\Xi\in C^1(U_\delta;\rline^n)$
such that \vspace{-1mm}
\begin{equation} \label{eq:xi_equil}
   f(\Xi(u),u) \m=\m 0 \FORALL u \in U_\delta, \vspace{-1mm}
\end{equation}
i.e., for each $u_0 \in U_\delta$, $\Xi(u_0)$ is an equilibrium point
that corresponds to the constant input $u_0$. Moreover, $\PPP$ is {\em
uniformly exponentially stable} around these equilibrium points. This
means that there exist $\e_0>0$, $\l>0$ and $m\geq 1$ such that for
each constant input function $u_0\in U_\delta$, the following holds
for the solutions of \rfb{eq:P}:
		
If $\norm{x(0)-\Xi(u_0)}\leq\e_0$, then for every $t\geq 0$,
\begin{equation} \label{eq:exp_stab}
   \norm{x(t)-\Xi(u_0)} \leq me^{-\l t}\norm{x(0) - \Xi(u_0)}.
\end{equation}
\end{assumption}
\end{mdframed}

\begin{remark} \label{rem:Jacob}
The uniform exponential stability condition above can be checked by
linearization: If the Jacobian matrices \vspace{-2mm}
\begin{equation} \label{eq:jac}
   A(u_0) \m=\m \left. \frac{\partial f(x,u)}{\partial x} \right|_
   {\nm\begin{array}{c} \scriptstyle x=\Xi(u_0)\vspace{-1mm}\\
   \scriptstyle u=u_0 \end{array}} \m\in\m \rline^{n\times n}
   \vspace{-1mm} \end{equation}
have eigenvalues bounded away from the right half-plane,
$$ \max\Re\sigma(A(u_0)) \m\leq\m \l_0 \m<\m 0 \FORALL u_0\in 
   U_\delta,$$
then $\PPP$ is uniformly exponentially stable, see (11.16) in
\cite{Khalil2002}. From the first part of Assumption \ref{ass_1},
$\max\Re\sigma(A(u_0))$ is a continuous function of $u_0$. Hence, if
this function is always negative, then by the compactness of
$U_\delta$, its maximum is also negative. Thus, for the uniform
exponential stability of $\PPP$ we only have to check that each of the
matrices $A(u_0)$ is stable (where $u_0\in U_\delta$).
\end{remark}

\begin{remark} \label{remark:C_2}
If there exists a function $\Xi$ with the properties \rfb{eq:xi_equil}
and \rfb{eq:exp_stab}, then automatically $\Xi\in C^2$ according to
the implicit function theorem, since $f\in C^2$. The stability
property \rfb{eq:exp_stab} guarantees that all the eigenvalues of the
Jacobians from \rfb{eq:jac} are non-zero, see Remark \ref{rem:Jacob}.
\end{remark}

\begin{mdframed}
\begin{assumption} \label{Ass2}
The system $\PPP$ satisfies Assumption \ref{ass_1} and, moreover, the
function
\begin{displaymath}\label{eq:G}
   G(u)\vcentcolon=g(\Xi(u)) \FORALL u\in U_\delta,
\end{displaymath}
is monotone increasing in the following sense: there exists $\mu>0$ 
such that
$$G(b)-G(a) \m\geq\m \mu(b-a) \FORALL a,b\in U_\delta,\ a<b \m.$$
\end{assumption}
\end{mdframed}

Note that the function $G\in{\rm Lip}(U_\delta,\rline)$. We denote
$y_{min}\vcentcolon= G(u_{min})$ and $y_{max}\vcentcolon=G(u_{max})$,
where clearly $y_{min}<y_{max}$. Moreover, using the notation
$Y=[y_{min}, y_{max}]$, for any $r\in Y$, we define
$u_r\vcentcolon=G^{-1}(r)$, which is well-defined in $U$ since $G$ is
continuous and strictly monotone.

\begin{remark}
If $\PPP$ is linear, described by the matrices $A,B,C$ in the usual
way ($\dot{x}=Ax+Bu$, $y=Cx$), then Assumption \ref{ass_1} reduces
to the fact that $A$ is stable. The functions $\Xi$ and $G$ from
Assumptions \ref{ass_1} and \ref{Ass2} are given by
$$ \Xi(u) \m=\m (-A)^{-1}Bu, \qquad G(u) \m=\m P(0)u,$$
where $P(s)=C(sI-A)^{-1}B$ is the plant transfer function. 
Assumption \ref{Ass2} is satisfied if and only if $P(0)=\mu>0$.
Note that the above conditions are the same 
conditions given in \cite{Morari1985} for the closed-loop stability of
a linear (stable) plant connected in feedback with a low-gain integral controller,
when $k>0$.
\end{remark}

Assumptions \ref{ass_1} and \ref{Ass2} are crucial for the stability
analysis presented in Section \ref{sec:stability_analysis}. In fact,
Assumption \ref{ass_1} guarantees the stability of the boundary-layer
(fast) system, while Assumption \ref{Ass2} guarantees the stability of
the reduced (slow) model. Note that $G$ is the steady-state
input-output map associated to $\PPP$.

In the following, we manipulate the closed-loop equations \rfb{Trump}
to rewrite them as a standard singular perturbation model (see
(11.47), (11.48) in \cite{Khalil2002}). Defining the function
\begin{equation} \label{eq:beta}
   \beta(x,u_I,k) \m\vcentcolon=\m f(x,u_I+\tau_p k(r-g(x)))
   -f(x,u_I),
\end{equation}
the closed-loop system \rfb{Trump} can be rewritten as
\begin{equation} \label{eq:closed_loop_system_beta}
   \begin{cases} \begin{gathered}
   \m \ \m\dot{x} \m=\m f(x,u_I)+\beta(x,u_I,k), \\ 
   \m\dot{u}_I \m=\m \Sscr(u_I,k(r-g(x))). \end{gathered} \end{cases}
\end{equation}
Note that $\beta(x,u_I,0)=0$ for all $(x,u_I)\in \mathbb{R}^{n+1}$. We define
\begin{equation} \label{eq:h}
h(x,u_I) \vcentcolon=\m \Sscr(u_I,r-g(x)),
\end{equation}
and then the closed-loop system
\rfb{eq:closed_loop_system_beta} can be rewritten as
\begin{equation} \label{eq:clos_h}
   \dot{u}_I \m=\m k\cdot h(x,u_I), \quad \dot{x} \m=\m f(x,u_I)
   + \beta(x,u_I,k).
\end{equation}
The functions $h(\cdot,\cdot)$ and $\beta(\cdot,\cdot,k)$ are defined
on $\rline^n\times U_\delta$. Using that $k>0$, we change the {\em
time-scale} of \rfb{eq:clos_h} introducing $s\vcentcolon=kt$ (it
is a ``slower'' time-scale because $k$ is small). Rewriting the system
\rfb{eq:clos_h} in the new time-scale $s$, we get
\begin{equation} \label{eq:clos_h_scale_s}
   \frac{\dd u_I}{\dd s} \m=\m h(x,u_I), \quad k\frac{\dd x}{\dd s} 
   \m=\m f(x,u_I)+\beta(x,u_I,k).
\end{equation}

To simplify the stability analysis, we move the equilibrium point of
the closed-loop system \rfb{eq:clos_h_scale_s} to the origin. It
follows from Assumption \ref{ass_1} and the notation introduced after
Assumption \ref{Ass2} that $(x_r,u_r)$ is an equilibrium point of
\rfb{eq:clos_h_scale_s}, where $x_r\vcentcolon =\Xi(u_r)$. Note that
$g(x_r)=r$. Introduce the variables
\begin{equation} \label{eq:x_tilde_u_tilde}
   \tilde{x}\m\vcentcolon=\m x-x_r, \qquad \tilde{u}_I
   \vcentcolon =\m u_I-u_r,
\end{equation}
and the functions
\begin{subequations}
   \begin{gather} \tilde{h}(\tilde{u}_I,\tilde{x})\vcentcolon=\m
   h(\tilde{x}+x_r,\tilde{u}_I+u_r), \label{eq:h_tilde} \\  
   \tilde{f}(\tilde{u}_I,\tilde{x})\vcentcolon=\m f(\tilde{x}+x_r,
   \tilde{u}_I+u_r), \label{eq:f_tilde} \\ \tilde{g}(\tilde{x})
   \vcentcolon=g(\tilde{x}+x_r), \label{eq:g_tilde} \\ 
   \tilde{\beta}(\tilde{u}_I,\tilde{x},k)\vcentcolon =\m \tilde{f}
   (\tilde{u}_I+\tau_p k(r-\tilde{g}(\tilde{x})),\tilde{x}) -
   \tilde{f}(\tilde{u}_I,\tilde{x}) \label{eq:b_tilde} \end{gather}
\end{subequations}
\vspace{-6mm}
$$=\m \beta(\tilde{x}+x_r,\tilde{u}_I+u_r,k).$$
The system \rfb{eq:clos_h_scale_s} can be rewritten as
\begin{equation} \label{f_h_tilde_s}
   \frac{\dd\tilde{u}_I}{\dd s} \m=\m \tilde{h}(\tilde{u}_I,
   \tilde{x}),\qquad k\frac{\dd\tilde{x}}{\dd s} \m=\m 
   \tilde{f}(\tilde{u}_I,\tilde{x})+\tilde{\beta}(\tilde{u}_I,
   \tilde{x},k),
\end{equation}
with an equilibrium point $(\tilde{u}_I,\tilde{x})=(0,0)$. For $k>0$
small, this is a \textit{standard singular perturbation model}
according to \cite[Sect.~11.5]{Khalil2002} (see equations (11.47)
and (11.48) there). Note that in our case the system 
\rfb{f_h_tilde_s} is autonomous. 

\begin{table}[t] \begin{center}
   \renewcommand{\arraystretch}{1.5}
   \caption{Correspondence of our notation with the one used in 
   \cite{Khalil2002}}\label{tb:corr}
   \begin{tabular}{ccccccccccc} \hline this paper & $\tilde{u}_I$ & 
   $\tilde{x}$ & $\tilde{h}$ & $\tilde{f}+\tilde{\beta}$  & 
   $k$ & $z$ & $\tilde{\Xi}$ & $s$ & $t$ & $\e_0$ \\ \hline
   \cite{Khalil2002} & $x$ & $z$ & $f$ & $g$ & $\e$ & $y$ & $h$ & $t$ & $\tau$ & 
   $r_0$\sbluff \\ \hline \end{tabular} \end{center} 
   \vspace{-4mm} \end{table}

Recalling $\Xi$ from Assumption 1, we introduce the function
\begin{equation} \label{eq:Xi_tilde}
   \tilde{\Xi}(\tilde{u}_I) \vcentcolon=\m \Xi(\tilde{u}_I+u_r)-x_r
\end{equation}
and following the guidelines of \cite{Khalil2002}, we define
\begin{equation} \label{eq:variable_change}
   z \vcentcolon =\m \tilde{x}-\tilde{\Xi}(\tilde{u}_I)=x-\Xi(u_I).
\end{equation}
Using the notation introduced above, we reformulate our system
\rfb{f_h_tilde_s} like (11.49), (11.50) of \cite{Khalil2002}, i.e.,
\begin{equation} \label{cl_loop_u_tilde}
   \frac{\dd\tilde{u}_I}{\dd s} \m=\m \tilde{h}(\tilde{u}_I,z +
   \tilde{\Xi}(\tilde{u}_I)),
\end{equation}\vspace{-6mm}
\begin{multline} \label{eq:cl_loop_z} 
   k\frac{\dd z}{\dd s} \m=\m\tilde{f}(\tilde{u}_I,z+\tilde{\Xi}
   (\tilde{u}_I)) + \tilde{\beta}(\tilde{u}_I,z+\tilde{\Xi}
   (\tilde{u}_I),k) \\ - k \frac{\dd\tilde{\Xi}}{\dd\tilde{u}_I}
   \tilde{h}(\tilde{u}_I,z+\tilde{\Xi}(\tilde{u}_I)),
\end{multline}
which has an equilibrium point at $(\tilde{u}_I,z)=(0,0)$. Finally, in
accordance with the change of variables \rfb{eq:x_tilde_u_tilde}, we
define
$$ \tilde{u}_{min} \vcentcolon=\m u_{min}-u_r \m,\qquad 
   \tilde{u}_{max} \vcentcolon=\m u_{max}-u_r$$ 
and the sets $\tilde{U}\vcentcolon=[\tilde{u}_{min},\tilde{u}_{max}]
\subset\rline$, and
$$ \tilde{U}_\delta \vcentcolon= [\tilde{u}_{min}-\delta,
   \tilde{u}_{max}+\delta] \m.$$
$\tilde{U}_\delta$ contains 0 in its interior. The state space of the 
closed-loop system \rfb{cl_loop_u_tilde}-\rfb{eq:cl_loop_z} is 
$\tilde{X}_\delta\vcentcolon=\tilde{U}_\delta\times\rline^n$. To
facilitate the comparison of our equations with those in
\cite{Khalil2002}, the relation between the notation in these two
works is shown in Table \ref{tb:corr}.

\begin{remark} \label{rem:stability_preserving}
Since $\tilde{\Xi}\in C^2(\tilde{U}_\delta,\rline^n)$, the change of
variables \rfb{eq:variable_change} is stability preserving, i.e., the
origin of \rfb{cl_loop_u_tilde}-\rfb{eq:cl_loop_z} is asymptotically
(exponentially) stable, if and only if the origin of \rfb{f_h_tilde_s}
is asymptotically (exponentially) stable.
\end{remark}

\begin{figure} 
\centering \includegraphics[width=0.3\textwidth]{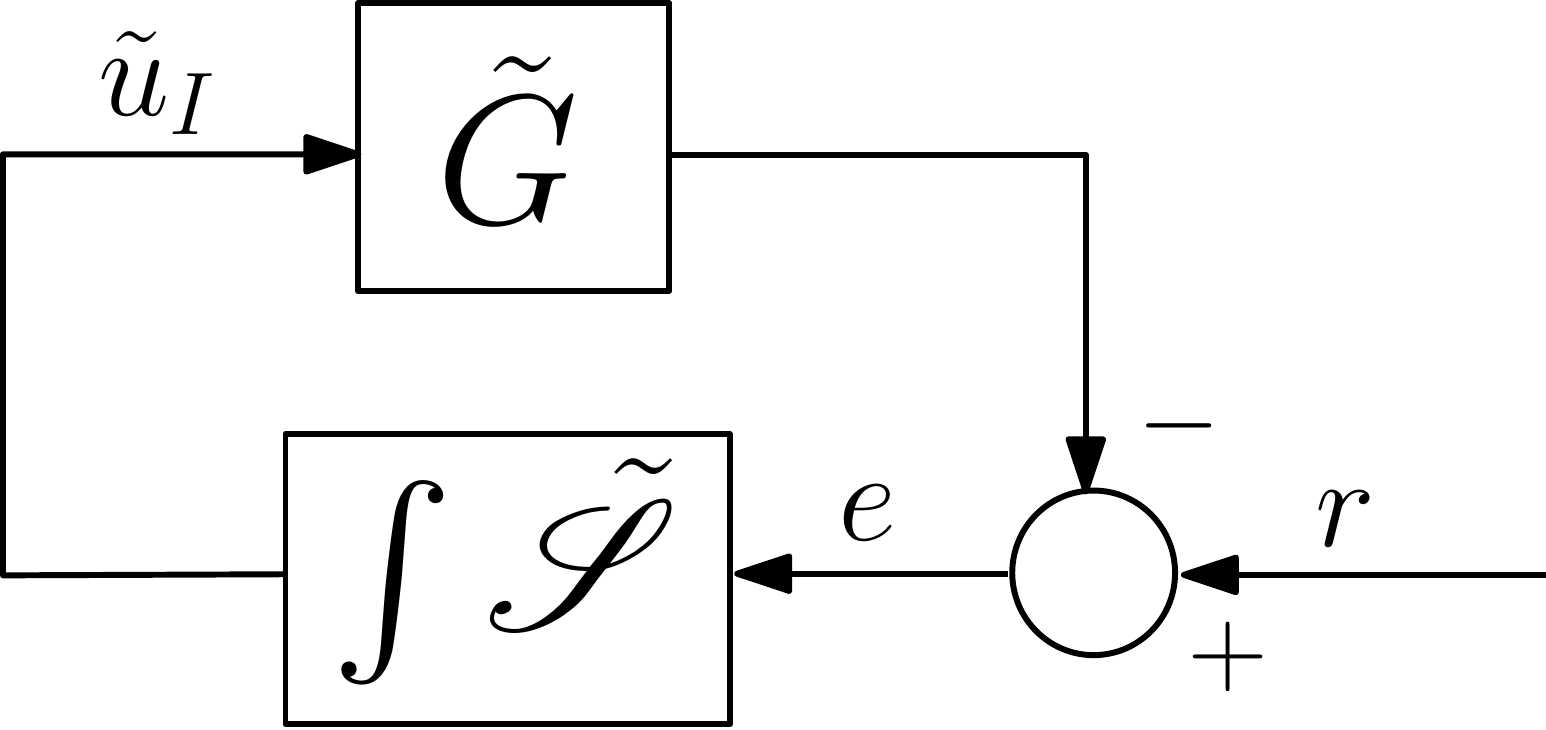}
\caption{Closed-loop representation of the reduced (slow) model 
\rfb{eq:reduced_model_final}.} \label{fig:red_mod}
\end{figure}

The \textit{reduced (slow) model} is obtained setting $k=0$ in
\rfb{eq:cl_loop_z}, and solving the resulting algebraic equation
$$0 \m=\m \tilde{f}(\tilde{u}_I,z+\tilde{\Xi}(\tilde{u}_I)),$$ for
which $z=0$ is a solution. Substituting $z=0$ in
\rfb{cl_loop_u_tilde}, the following reduced (slow) model is obtained:
\begin{equation} \label{eq:reduced_model}
   \frac{\dd\tilde{u}_I}{\dd s} \m=\m \tilde{h}(\tilde{u}_I,
   \tilde{\Xi}(\tilde{u}_I)).
\end{equation}
Recall $\Sscr$ from \rfb{eq:sat} and $G$ from Assumption \ref{Ass2}.
Defining
\begin{equation*} \label{eq:S_G_tilde}
   \begin{gathered} \tilde{\mathscr{S}}(\tilde{u}_I,\cdot)
   \vcentcolon=\m \Sscr(\tilde{u}_I+u_r,\cdot), \quad \tilde{G}
   (\tilde{u}_I)\vcentcolon=\m G(\tilde{u}_I+u_r),\end{gathered}
\end{equation*}
the reduced model \rfb{eq:reduced_model} can be described equivalently
(recall the definitions of $\tilde{h}$ from \rfb{eq:h_tilde} and of
$h$ from \rfb{eq:h}) by
\begin{equation} \label{eq:reduced_model_final}
   \frac{\dd\tilde{u}_I}{\dd s} \m=\m \tilde{\Sscr}(\tilde{u}_I,r-
   \tilde{G}(\tilde{u}_I)).
\end{equation}
The reduced closed-loop system is shown in Fig. \ref{fig:red_mod}. 
Note that
\begin{equation*} \label{eq:h_tilde_S_tilde}
   \tilde{h}(\tilde{u}_I,\tilde{x}) \m=\m \tilde{\Sscr}(\tilde{u}_I,
   r-\tilde{g}(\tilde{x})).
\end{equation*}

The \textit{boundary-layer (fast) system} is obtained by rewriting the
second equation in \rfb{f_h_tilde_s} in the original fast time-scale
$t$ as
\begin{equation*} \label{eq:fast_system}
   \dot{z} \m=\m \tilde{f}(\tilde{u}_I,z+\tilde{\Xi}(\tilde{u}_I)) +
   \tilde{\beta}(\tilde{u}_I,z+\tilde{\Xi}(\tilde{u}_I),k).
\end{equation*} 
Taking here $k=0$, we obtain the boundary-layer system
\begin{equation} \label{eq:fast_system_k_0}
   \dot{z}\m=\m \tilde{f}(\tilde{u}_I,z+\tilde{\Xi}(\tilde{u}_I)),
\end{equation}
where $\tilde{u}_I\in\tilde{U}_\delta$ is treated as a fixed parameter.

\section{Stability analysis via singular perturbations}
         \label{sec:stability_analysis} 

To perform the stability analysis of the closed-loop system
\rfb{f_h_tilde_s}, we would like to use Theorem 11.4 from
\cite{Khalil2002}. However, the assumptions of this theorem are
not met by our model. Therefore, we state a slightly modified version
of it, which holds in our framework, together with a
proof. Note that a similar result, with asymptotic stability in place
of exponential one, is provided in Theorem 5.1 from
\cite[Chapt.~7]{Kokotovic1999}.

The analysis will be performed by first finding two distinct Lyapunov
functions $V$ and $W$, respectively for the reduced (slow) model
\rfb{eq:reduced_model} and the boundary-layer (fast) system
\rfb{eq:fast_system_k_0}. Then, our modified version of Theorem 11.4
from \cite{Khalil2002} will be stated, see Theorem \ref{thm:11.4}. The
stability of the closed-loop system \rfb{Trump}, with the region of
attraction containing the curve $\{(\Xi(u),u)\m|\m u\in U\}$, will be
presented in Theorem \ref{thm:stab}.

\subsection{Stability of the reduced (slow) model}
\label{subsec:5A} 

Our Theorem \ref{thm:11.4} requires the existence of a Lyapunov
function $V$ for \rfb{eq:reduced_model} (defined on 
$\tilde{U}_\delta$) such that
\begin{equation} \label{eq:V_condit}
   \begin{gathered} c_1\tilde{u}_I^2\leq V(\tilde{u}_I)\leq c_2
   \tilde{u}_I^2, \\ \frac{\dd V}{\dd \tilde{u}_I}\tilde{h}
   (\tilde{u}_I,\tilde{\Xi}(\tilde{u}_I))\leq-c_3\abs{\tilde{u}_I}^2,
   \ \ \ \bigg|\frac{\dd V}{\dd \tilde{u}_I}\bigg|\leq c_4
   \abs{\tilde{u}_I},\end{gathered}
\end{equation}
for all $\tilde{u}_I\in\tilde{U}_\delta$, where $c_1,\dots c_4$ are 
positive constants. Consider the candidate
\begin{equation} \label{eq:lyap_reduced_model_1}
   V(\tilde{u}_I) \m=\m \frac{1}{2}\tilde{u}_I^2 \FORALL \tilde{u}_I
   \in\tilde{U}_\delta.
\end{equation}
It follows from \rfb{eq:reduced_model_final} that 
$$ \frac{\dd V}{\dd s} \m=\m \tilde{u}_I\frac{\dd\tilde{u}_I}
   {\dd s} \m=\m \tilde{u}_I\tilde{\mathscr{S}}(\tilde{u}_I,r-
   \tilde{G}(\tilde{u}_I)).$$
If $\tilde{u}_I(s)=0$, then the system \rfb{eq:reduced_model_final}
has reached its equilibrium point (recall that $r=\tilde{G}(0)$). If
$\tilde{u}_I(s)>0$, then recalling that $\tilde{G}$ is monotone 
increasing (see Assumption \ref{Ass2}), it follows that
\begin{displaymath}
   \frac{\dd\tilde{u}_I}{\dd s}=\tilde{G}(0)-\tilde{G}
   (\tilde{u}_I(s))<0,
\end{displaymath}
therefore
\begin{displaymath}
   \frac{\dd V}{\dd s} \m=\m \tilde{u}_I\frac{\dd\tilde{u}_I}
   {\dd s} \m<\m 0 \m.
\end{displaymath}
If $\tilde{u}_I(s)<0$, then by a similar argument the same conclusion
is obtained. This proves that \rfb{eq:lyap_reduced_model_1} is indeed
a Lyapunov function for \rfb{eq:reduced_model}, defined for all
$\tilde{u}_I\in\tilde{U}_\delta$.

\begin{remark}
Note that the above argument tells us also that if $\tilde{u}_I(0)
\in\tilde{U}$, then the saturating integrator $\int\tilde
{\mathscr{S}}$ of the reduced model in Fig. \ref{fig:red_mod} is
behaving like a standard integrator.
\end{remark}

From Assumption \ref{Ass2}, $\abs{\tilde{G}(\tilde{u}_I)-\tilde{G}
(0)}\geq\mu\abs{\tilde{u}_I}$, for all $\tilde{u}_I\in\tilde{U}
_\delta$. Then we can write
\begin{displaymath}
   \frac{\dd V}{\dd s} = -\left| \tilde{u}_I\sbluff \left[\tilde{G}(0)
   -\tilde{G}(\tilde{u}_I)\right] \right| \leq -\abs{\tilde{u}_I}\mu
   \abs{\tilde{u}_I}= -\mu\abs{\tilde{u}_I}^2,
\end{displaymath}
for all $\tilde{u}_I\in\tilde{U}_\delta$. Therefore, the conditions 
\rfb{eq:V_condit} holds with
\begin{displaymath}
   \begin{gathered} c_1=c_2=\half, \quad c_3=\mu \quad \text{and} 
   \quad c_4=1. \end{gathered}
\end{displaymath}

\subsection{Stability of the boundary-layer (fast) system}
\label{subsec:5B} 

Our Theorem \ref{thm:11.4} requires the existence of a Lyapunov
function $W$ for \rfb{eq:fast_system_k_0} (defined on $\tilde{U}
_\delta\times B_{\e_0}$) such that
\begin{equation} \label{eq:W_lemma_9_8}
   \begin{gathered} b_1\norm{z}^2\leq W(\tilde{u}_I,z)\leq b_2
   \norm{z}^2, \\ \frac{\partial W}{\partial z} \tilde{f}
   (\tilde{u}_I,z+\tilde{\Xi}(\tilde{u}_I)) \leq -b_3\norm{z}^2, \\
   \bigg\|\frac{\partial W}{\partial z}\bigg\| \leq b_4\norm{z},\quad
   \bigg\|\frac{\partial W}{\partial \tilde{u}_I}\bigg\| \leq b_5
   \norm{z}^2,\end{gathered}
\end{equation}
for all $(\tilde{u}_I,z)\in \tilde{U}_\delta\times B_{\e_0}$ (recall
$\e_0$ from Assumption \ref{ass_1}), where $b_1,\dots b_5$ are
positive constants and $B_{\e_0}$ denotes the closed ball of radius
$\e_0$ in $\rline^n$. We use Lemma 9.8 of \cite{Khalil2002}, which,
under Assumption \ref{ass_1} and some smoothness requirements on
$\tilde{f}$, guarantees the existence of $W$ such that the conditions
in \rfb{eq:W_lemma_9_8} hold. The aforementioned requirements on
$\tilde{f}$ are: Denoting $p(z,\tilde{u}_I)\vcentcolon=\tilde{f}
(\tilde{u}_I,z+\tilde{\Xi}(\tilde{u}_I))$ (this function is denoted by
$g$ in Lemma 9.8 of \cite{Khalil2002}), it should hold that
\begin{equation} \label{eq:lemma_9_8_jac}
   \bigg\|\frac{\partial p}{\partial z}(z,\tilde{u}_I)\bigg\|
   \leq L_1 \quad \text{and} \quad \bigg\|\frac{\partial p}
   {\partial \tilde{u}_I}(z,\tilde{u}_I)\bigg\|\leq L_2\|z\|,
\end{equation}
for all $(z,\tilde{u}_I)\in B_\rho\times\tilde{U}_\delta$ (where
$\e_0\leq\rho<\infty$). The first condition follows from $p\in C^2$
and $B_\rho\times\tilde{U}_\delta$ being compact. For the second
condition, for any $j\in\{1,2,\dots n\}$, we define $F_j(z,\tilde{u}
_I)\vcentcolon=\frac{\partial p_j}{\partial\tilde{u}_I}(z,\tilde{u}
_I)\in\rline$. Since $\tilde{f}(\tilde{u}_I,\tilde{\Xi}(\tilde{u}_I))
=0$ for all $\tilde{u}_I\in\tilde{U}_\delta$, then $F_j(0,\tilde{u}
_I)=0$ for all $\tilde{u}_I\in\tilde{U}_\delta$. For every fixed 
$(z,\tilde{u}_I)\in B_{r}\times\tilde{U}_\delta$, we introduce the
function $\tilde{F}_j:[0,1]\to\rline$ such that $\tilde{F}_j(\sigma)=
F_j(\sigma z,\tilde{u}_I)$, hence $\tilde{F}_j(0)=0$ and $\tilde{F}_j
(1)=F_j(z,\tilde{u}_I)$. According to the mean value theorem, there
exists $\xi\in(0,1)$ such that $F_j(z,\tilde{u}_I)=\tilde{F}_j(1)-
\tilde{F}_j(0)=\tilde{F}_j'(\xi)$. Therefore, we get that for all
$(z,\tilde{u}_I)\in B_r\times\tilde{U}_\delta$,
$$ \abs{F_j(z,\tilde{u}_I)}\leq\bigg\|\frac{\partial F_j}
   {\partial z}(\xi z,\tilde{u}_I)\bigg\| \cdot\norm{z}=\bigg\|
   \frac{\partial^2 p_j}{\partial z \partial
   \tilde{u}_I}(\xi z,\tilde{u}_I)\bigg\|\cdot\norm{z},$$
which holds because $p\in C^2$ and $B_{r}\times\tilde{U}_\delta$ is
compact. Recall that $F_j(z,\tilde{u}_I)$ is the $j$-th component of
$\frac{\partial p}{\partial \tilde{u}_I}(z,\tilde{u}_I)$, then the
above inequality implies the second estimate in 
\rfb{eq:lemma_9_8_jac}. Therefore, we can apply Lemma 9.8 of
\cite{Khalil2002}, which yields the existence of $W$ and $b_1\dots 
b_5$ such that \rfb{eq:W_lemma_9_8} holds.

\vspace{-3mm}
\subsection{Closed-loop stability analysis}
            \label{subsec:4p4} 

First we state our modified version of Theorem 11.4 from
\cite{Khalil2002}. We have strived to formulate a self-contained
theorem. For this, we forget about the specific meaning of the
notation $\tilde{f},\tilde{\beta},\tilde{h},\tilde{\Xi},
\tilde{U}_\delta$ introduced earlier in this paper and we introduce
these functions in our Theorem \ref{thm:11.4}, while stating all the
properties that are needed for this theorem.

\smallskip
{\color{blue}
\begin{theorem} \label{thm:11.4}
Let $\tilde{U}_\delta$ be a finite closed interval containing 0 in its
interior. Consider the singularly perturbed system \rfb{f_h_tilde_s},
with state space $\tilde{X}_\delta=\tilde{U}_\delta\times\rline^n$,
where $\tilde{f}\in C^2(\tilde{X}_\delta; \rline^n)$, $\tilde{\beta}
\in {\rm Lip}(\tilde{X}_\delta\times\rline; \rline^n)$, and $\tilde{h}:
\tilde{X}_\delta\to\rline$ is uniformly Lipschitz in the second
argument. We assume that \rfb{f_h_tilde_s} has a unique local solution
for every initial state in $\tilde{X}_\delta$ and for every $k>0$.
There is also a function $\tilde{\Xi}\in C^1(\tilde{U}_\delta;
\rline^n)$.  Assume the following:
\begin{itemize}
\item[a.] $\tilde{f}(0,0)=0$ and $\tilde{\beta}(0,0,k)=0$ for all
   $k>0$.
\item[b.] For every $\tilde{u}_I\in\tilde{U}_\delta$, the equation
   $\tilde{f}(\tilde{u}_I,\tilde{x})=0$ has an isolated root
   $\tilde{x}=\tilde{\Xi}(\tilde{u}_I)$ and $\tilde{\Xi}(0)=0$.
\item[c.] $\tilde{\beta}(\tilde{u}_I,\tilde{x},0)=0$ for all 
   $(\tilde{u}_I,\tilde{x})\in\tilde{X}_\delta$, and $\tilde{h}
   (\tilde{u}_I,0)=0$ for all $\tilde{u}_I\in\tilde{U}_\delta$.
\item[d.] There exists a Lyapunov function $V$ for the reduced system
   \rfb{eq:reduced_model} such that \rfb{eq:V_condit} holds.
\item[e.] There exists a Lyapunov function $W$ for the boundary-layer
   system \rfb{eq:fast_system_k_0} such that \rfb{eq:W_lemma_9_8} 
   holds.
\end{itemize}

Then there exists $\kappa>0$ such that for all $k\in(0,\kappa)$, the
origin of \rfb{f_h_tilde_s} is exponentially stable. Moreover, the
function
\begin{equation} \label{eq:nu_Lyap}
   \nu(\tilde{u}_I,z)\vcentcolon=V(\tilde{u}_I)+W(\tilde{u}_I,z)
\end{equation}
is a Lyapunov function for the system \rfb{cl_loop_u_tilde}-%
\rfb{eq:cl_loop_z}, obtained from \rfb{f_h_tilde_s} by the change of
variables \rfb{eq:variable_change}.
\end{theorem}}

\smallskip
For the proof see Appendix \ref{app:B}.

We now return to the context of our problem introduced in Section
\ref{sec2}, and we use all the notation introduced in Sections
\ref{sec2} and \ref{sec3}. We use the above ``general'' theorem to
prove that the equilibrium point $(\Xi(u_r),u_r)$ of the closed-loop
system \rfb{Trump} is locally exponentially stable, with its region of
attraction containing the curve $\{(\Xi(u),u)\m|\m u\in U\}$. If
the initial state is in the region of attraction, then the output $y$
of $\PPP$ converges to $r\in Y$ (recall that $Y=[y_{min},y_{max}]$ 
and $u_r=G^{-1}(r)$).

\smallskip
{\color{blue}
\begin{theorem} \label{thm:stab}
Consider the closed-loop system \rfb{Trump}, where $\PPP$ satisfies
Assumption \ref{Ass2}. Then there exists a $\kappa>0$ such that if the
gain $k\in(0,\kappa]$, then for any $r\in Y$, $(\Xi(u_r),u_r)$ is a
(locally) exponentially stable equilibrium point of the closed-loop
system \rfb{Trump}, with state space $X=\rline^n\times U$.

If the initial state $(x(0),u_I(0))\in X$ of the
closed-loop system satisfies $\norm{x(0)-\Xi(u_I(0))}\leq\e_0$, then
$$ x(t)\rarrow\Xi(u_r), \qquad u(t)\rarrow u_r, \qquad y(t)
   \rarrow r,$$
and this convergence is at an exponential rate. 
\end{theorem}}

\smallskip
\begin{IEEEproof} Take $r\in Y$. Consider the singularly perturbed
system \rfb{f_h_tilde_s} with state space $\tilde{X}_\delta=\tilde{U}
_\delta\times\rline^n$. Recall from Subsections \ref{subsec:5A} and
\ref{subsec:5B} that (thanks to Assumption \ref{Ass2}) there are
functions $V$ and $W$ satisfying the conditions \rfb{eq:V_condit} and
\rfb{eq:W_lemma_9_8}. According to Theorem \ref{thm:11.4} there exists
a $\kappa_r>0$ such that the origin of \rfb{f_h_tilde_s} is (locally)
exponentially stable for all $k<\kappa_r$. Moreover, $\nu$ from
\rfb{eq:nu_Lyap} is a Lyapunov function for the closed-loop system
\rfb{cl_loop_u_tilde}-\rfb{eq:cl_loop_z}. From
\rfb{eq:lyap_reduced_model_1}, it is clear that
\begin{displaymath}
   V(\tilde{u}_I)\leq\frac{1}{2}\max\{(\tilde{u}_{min}-\delta)^2,
   (\tilde{u}_{max}+\delta)^2\}=\vcentcolon u_{M},
\end{displaymath}
for all $\tilde{u}_I\in\tilde{U}$, and from \rfb{eq:W_lemma_9_8}, it
follows that
\begin{equation} \label{eq:W_bound}
   W(\tilde{u}_I,z)\leq b_2\e_0^2 \FORALL (\tilde{u}_I,z)\in 
   \tilde{U}_\delta\times B_{\e_0}.
\end{equation}
We introduce the compact positively-invariant set \vspace{-2mm}
\begin{multline} \label{eq:L}
   L\vcentcolon =\m \bigg\{(\tilde{u}_I,z) \in \tilde{U}_\delta
   \times\rline^n\ |\ \nu(\tilde{u}_I,z) \leq u_{M}+b_2\e_0^2\bigg\}.
\end{multline}
From the above reasoning $\tilde{U}_\delta\times B_{\e_0}\subset L$,
and therefore, $\tilde{U}_\delta\times B_{\e_0}$ is contained in the
domain of attraction of the origin of 
\rfb{cl_loop_u_tilde}-\rfb{eq:cl_loop_z} for all $k<\kappa_r$.
Finally, the change of variables \rfb{eq:variable_change} is stability
preserving (see Remark \ref{rem:stability_preserving}), hence the
equilibrium point $(\Xi(u_r),u_r)$ of \rfb{Trump} is locally
exponentially stable with region of attraction containing all
$(x_0,u_0)\in\rline^n\times U_\delta$ such that $\|x_0-\Xi(u_0)\|\leq
\e_0$. Since $y(t)=g(x(t))$ and $g$ is locally Lipschitz, we have that 
$y(t)$ converges to $g(\Xi(u_r))= G(u_r)=r$ at an exponential rate.

We now want to find a $\kappa>0$ independent from the specific $r\in
Y$. Substituting $b_3$ in \rfb{eq:W_lemma_9_8} with $\tilde{b}_3=
\frac{b_3}{2}$, there exists an $\e_r>0$ such that
\rfb{eq:W_lemma_9_8}, with $W$ from above and $\tilde{b}_3$ in place
of $b_3$, still holds for all $r_1\in(r-\e_r,r+\e_r)=\vcentcolon
Y_r$. Indeed, the first and third lines of \rfb{eq:W_lemma_9_8} are
independent of $r$, while the second one is true in a neighborhood of
$r$ because $\tilde{f}$ and $\tilde{\Xi}$ are continuous functions of
$r$ (see \rfb{eq:f_tilde}, \rfb{eq:Xi_tilde}). Similarly, the
stability of the reduced system \rfb{eq:reduced_model} is guaranteed
for all $r_1\in Y_r$. Hence we can apply Theorem \ref{thm:11.4} to
find a $\kappa_r>0$ such that $(\Xi(u_r),u_r)$ is an exponentially
stable equilibrium point of the closed-loop system \rfb{Trump}, for
all $r_1\in Y_r$. Since $Y$ is compact and the sets \{$Y_r\ | \ r\in
Y$\} are an open covering of $Y$, we can extract a finite cover
\{$Y_r\ | \ r\in Y_F$\}. Therefore, it is enough to take $\kappa=\min
\{\kappa_r\ |\ r\in Y_F\}$.
\end{IEEEproof}

\begin{remark}\label{rmk:bound_kappa} In Theorem
\ref{thm:stab} we have proved the existence of an upper bound
$\kappa>0$ for the controller gain $k$, guaranteeing the (exponential)
stability of the closed-loop system equilibrium points for varying
$r\in Y$. For the interested readers, we mention that a formula for
$\kappa$ can be found following the steps of the proof of Theorem
\ref{thm:11.4} (see Appendix \ref{app:B}). Indeed, for the inequality
(obtained before \rfb{eq:exp_bound} in Appendix \ref{app:B})
\begin{displaymath} 
   \frac{\dd\nu}{\dd s} \m\leq\m -\bbm{ \abs{\tilde{u}_I} \\ \norm{z}
   }^\top \bbm{ c_3 & -d_6 \\ -d_6 & \frac{b_3}{k}-d_5} \bbm{
   \abs{\tilde{u}_I} \\ \norm{z}} \m<\m 0
\end{displaymath} 
to hold, we need $k<\frac{c_3b_3}{d_6^2+c_3d_5}=\vcentcolon\kappa_r$.
From here, following the steps of the last part of the proof of
Theorem \ref{thm:stab}, we can obtain a $\kappa>0$ independent from
the specific $r\in Y$. A similar procedure can be found in other works
using singular perturbation methods, see, for instance,
\cite[Theorem~4.2]{Lorenzetti2020} (for the case $\tau_p=0$),
\cite[Lemma~3.2]{Desoer1985}, \cite[Theorems~11.3,~11.4]{Khalil2002},
\cite[Theorem~2.1 (Chapt.~7)]{Kokotovic1999},
\cite[Theorem~3.1]{Simpson-Porco2020}.

Computing $\kappa$ using the method sketched above will lead to a very
complicated and extremely conservative expression, which is not of
practical help for tuning a real controller (an unnecessarily small
controller gain would cause a slow acting feedback and the resulting
closed-loop system would exhibit a large settling time). For this
reason, we have omitted the details of such method and, in practice,
the tuning should be performed through simulation experiments.  (The
excessive conservativeness of the upper bound $\kappa$ is a common
feature of all the stability analysis methods based on singular
perturbations theory, see the references mentioned above.)
\end{remark}

We show that if $r$ is in a class of step functions, then the
stability result of the above theorem still holds.

{\color{blue}
\begin{proposition} \label{prop:T}
Consider the closed-loop system \rfb{Trump}, where $\PPP$ satisfies
Assumption \ref{Ass2} and $k\in(0,\kappa)$ as in Theorem
\ref{thm:stab}. Then there exists a $\bar{T}>0$ with the following
property: Let $r$ be a step function with values in $Y$ with finitely
many discontinuity (jump) points. If the time difference between any
two jump points is at least $\bar{T}$, then the statement in the last
sentence of Theorem \ref{thm:stab} remains true.
\end{proposition}}

\begin{IEEEproof}
Recall the estimate \rfb{eq:exp_bound} from Appendix \ref{app:B}.
By rewriting it only for $z$, we get
\begin{equation}\label{eq:exp_b}
\norm{z(s)}\leq K_1e^{-\gamma s}\norm{z(0)}\leq K_1e^{-\gamma s}\e_0.
\end{equation}
We are interested in computing the time $\bar{T}$ needed by the state
trajectory $z(t)$ to enter again in $B_{\e_0}$, after the transient
due to the boundary-layer (fast) system trajectory of
\rfb{eq:fast_system_k_0}. From \rfb{eq:exp_b}, it is easy to obtained a
lower bound on $s$ as \vspace{-1mm}
\begin{displaymath} \label{eq:s_bound}
   s \m\geq\m \frac{1}{\gamma}\log K_1,
\end{displaymath}
and, recalling that $s=kt$, we get a lower bound on $t$ as 
\vspace{-1mm}
$$ t \m\geq\m\bar{T}\vcentcolon=\frac{1}{\gamma k}\log K_1.$$
We have that $\norm{z(t)}\leq\e_0$ for all $t\geq\bar{T}$. Now we 
reset the time to be $0$ at the first discontinuity point of $r$,
and then the same argument applies again and again. When we reach 
$t_1$, the last discontinuity point of $r$, then Theorem 
\ref{thm:stab} can be applied with initial state $(x(t_1),u_I(t_1))\in
X$ (satisfying $\|x(t_1)-\Xi(u_I(t_1))\|=\|z(t_1)\|\leq\e_0$) and with
the constant reference $r(t_1)\in Y$.
\end{IEEEproof}

\section{Global asymptotic stability} \label{sec:5} 

In the previous sections, we have derived local results for the
closed-loop system \rfb{Trump} shown in Fig. \ref{fig:clos}: the
existence and uniqueness of solutions has been stated in Proposition
\ref{prop:uniqueness}, while constant reference tracking has been
proved in Theorem \ref{thm:stab}. At this point, a legitimate question
is which assumptions are required for these results to hold in a
global framework. In our view, the concept of \textit{input-to-state
stability (ISS)} (see \cite{Sontag2008,Sontag1996}) enables such an
extension in a natural way. In fact, we will show how an additional
ISS type requirement on the open-loop plant $\PPP$ is enough to obtain
global stability results for the closed-loop system \rfb{Trump}.

The following assumption requires that our plant $\PPP$ from
\rfb{eq:P} satisfies the \textit{asymptotic gain (AG)} property
introduced in \cite{Sontag1996}, around each equilibrium point
$\Xi(u_0)$, with $u_0\in U$.

\textbf{Notation.} Let $\rline_+=[0,\infty)$ and denote\vspace{-1mm}
$$ \Kscr_\infty\vcentcolon= \left\{\gamma\in C(\rline_+;\rline_+)\m
   \bigg|\ \gamma(0)=0,\begin{array}{c} \gamma\text{ increasing}, \\
   \gamma\text{ unbounded} \end{array}\nm\right\},$$ \vspace{-1mm}
$$ \Kscr\nm\Lscr\vcentcolon=\nm \left\{\alpha\in C(\rline_+\nm\times
   \rline_+;\rline_+) \m\bigg| \begin{array}{c} \alpha(\cdot,s)
   \in\Kscr_\infty \m\forall s\geq 0, \\ \alpha(r,\cdot) 
   \text{ is decreasing} , \\ {\lim_{s\rarrow\infty}\alpha(r,s) = 0}
   \end{array} \nm\right\}.$$\vspace{-1mm}

\begin{mdframed}
\begin{assumption} \label{ass_3}
The plant $\PPP$ satisfies Assumption \ref{ass_1} and there exists
$\gamma\in\Kscr_\infty$ such that for each $x_0\in\rline^n$, each
$u_0\in U_\delta$ and all $u\in C(\rline_+;\rline)$, the following
holds: Denoting the solution of \rfb{eq:P} corresponding to the
initial state $x_0$ and the input $u$ by $x$, we have
$$ \limsup_{t\rarrow\infty}\norm{x(t)-\Xi(u_0)} \leq\m \gamma
   \left(\limsup_{t\rarrow\infty}\abs{u(t)-u_0}\right).$$
\end{assumption} \end{mdframed}

{\color{blue}
\begin{proposition}\label{prop:ISS}
If $\PPP$ satisfies Assumption \ref{ass_3}, then for each $u_0\in
U_\delta$ there exists $\alpha_{u_0}\in\Kscr\nm\Lscr$ such that for
all $x_0\in\rline^n$ and all $u\in C(\rline_+,\rline)$, the following
holds: Denoting the solution of \rfb{eq:P} corresponding to the
initial state $x_0\in\rline^n$ and the continuous input $u$ by $x$,
\begin{multline} \label{eq:ISS}
   \norm{x(t)-\Xi(u_0)} \m\leq\m \alpha_{u_0}\left(\norm{x_0-\Xi
   (u_0)},t\right)\\ +\gamma\left(\sup_{\sigma\in[0,t]}
   \norm{u(\sigma)-u_0}\right) \ \ \forall\m t\geq 0.
\end{multline}
\end{proposition}}

Note that \rfb{eq:ISS} means that $\PPP$ is ISS around each
equilibrium point corresponding to $u_0\in U_\delta$,
uniformly with respect to $u_0$.

\begin{IEEEproof}
Take $u_0\in U_\delta$. From Assumption \ref{ass_1}, it follows that
$\PPP$ is \textit{0-stable} (see \cite{Sontag1996}), around the
equilibrium point $\Xi(u_0)$. (Note that in \cite{Sontag1996} the
0-stability property is defined in terms of (Lyapunov) stability of
the origin for zero input, while, in our framework, it refers to the
(Lyapunov) stability of the equilibrium point $\Xi(u_0)$ for $\PPP$
subjected to constant input $u_0$). On the other hand, Assumption
\ref{ass_3} guarantees that $\PPP$ satisfies the AG property.
Therefore, to complete the proof, is enough to apply Theorem 2 from
\cite{Sontag2008}, which states that a systems is ISS if and only if
it is 0-stable and AG.
\end{IEEEproof}

Besides Assumption \ref{ass_3}, for the results of this section, we
need that the signal $u$ in the closed-loop system from Fig.
\ref{fig:clos} is constrained to a compact set. To guarantee this, we
assume that either $g$ is bounded or $\tau_p=0$. Note that the
boundedness of $g$ is true for most physical systems.

{\color{blue}
\begin{proposition}\label{prop:ISS_sol}
Assume that $\PPP$ satisfies Assumption \ref{ass_3},
and at least one of the following conditions: $g$ is bounded or
$\tau_p=0$. Take $k>0$ and $r\in Y$. Then the closed-loop system
\rfb{Trump} has a unique state trajectory on the interval
$[0,\infty)$, for any initial state in $X_\delta=\rline^n \times
U_\delta$. Moreover, at any time $t\geq 0$, the state $(x(t), u_I(t))$
depends continuously on the initial state $(x(0),u_I(0))\in X$.
\end{proposition}}

For the proof see Appendix \ref{app:C}.

We now state a global stability result for the closed-loop system
\rfb{Trump}, extending the local one from Theorem \ref{thm:stab}.

{\color{blue}
\begin{theorem} \label{thm:ISS_stability}
We work under the assumptions of Proposition \ref{prop:ISS_sol} and, additionally,
we assume that ${\mathbf P}$ satisfies Assumption \ref{Ass2}. Then
there exists a $\kappa^*>0$ such that for any $k\in(0,\kappa^*]$, for
any $r\in Y$, and for any initial state $(x_0,u_0)\in X$, the state
trajectory $(x,u_I)$ of the closed-loop system \rfb{Trump} satisfies
$$ x(t)\rarrow\Xi(u_r),\qquad u(t)\rarrow u_r,\qquad y(t)\rarrow r,$$
and this convergence is at an exponential rate.
\end{theorem}}

For the proof see Appendix \ref{app:D}.

Note that the above exponential convergence holds uniformly on any
compact set of initial states, but not globally (see Appendix 
\ref{app:D} for more details on the choice of $\kappa^*$).

\section{Examples} \label{sec:ex} 

In this section, two examples are presented. In Example \ref{ex:1}, we
illustrate our main global results from Section \ref{sec:5} for the
closed-loop system formed by a simple nonlinear plant $\PPP$ and a PI
controller with a saturating integrator, connected as in Fig.
\ref{fig:clos}. In Example \ref{ex:2}, we address a problem of
practical relevance: the control of the local voltage for an inverter
connected to the power grid.  In particular, we show how the PI
controller with a saturating integrator improves the recovery of the
system after a grid fault. Due to the high complexity of this system,
some details will be omitted.

\begin{figure} 
   \centering \includegraphics[width=0.45\textwidth]{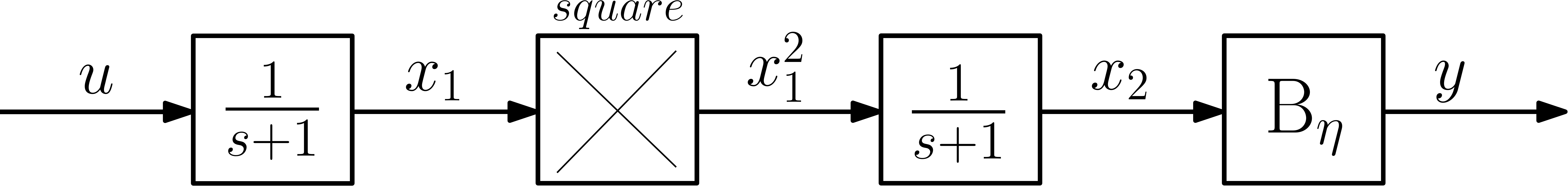}
   \caption{Equivalent block diagram of \rfb{eq:P1}.} \label{fig:P_1}
\end{figure}

\subsection{Control of a simple nonlinear plant} \label{ex:1} 

Consider he plant $\PPP$ shown in Fig. \ref{fig:P_1}, described by
\begin{equation} \label{eq:P1}
   \dot{x}_1 \m=\m -x_1+u, \quad \dot{x}_2 \m=\m x_1^2-x_2, \quad 
   y \m=\m {\mathrm B}_\eta(x_2),
\end{equation}
where $x\vcentcolon=\sbm{x_1\\ x_2}\in\rline^2$ is the state and
${\mathrm B}_\eta$ is the usual saturation function, defined for 
$\eta>0$ by
\begin{equation} \label{eq:B_eta}
   \mathrm{B}_\eta(z) \m=\m \begin{cases} -\eta & \text{if} \quad z
   \leq -\eta, \\ \ \ z & \text{if} \quad z \in (-\eta,\eta), \\
   \ \ \eta & \text{if} \quad z\geq \eta. \end{cases}
\end{equation}

The functions $\Xi$ and $G$ from Assumptions \ref{ass_1}
and \ref{Ass2} are
\begin{equation*} \label{eq:Xi_G}
   \Xi(u) \m=\m \begin{bmatrix} u \\ u^2
   \end{bmatrix} \quad \text{and} \quad G(u)=u^2,
\end{equation*}
for all $u\in U_\delta$, where $U_\delta\vcentcolon=[u_{min}-\delta,u_{max}+
\delta]$. We take positive numbers $u_{min}<u_{max}$ and $\delta$ satisfying
$$0 \m<\m u_{min}-\delta \m,\quad (u_{max}+\delta)^2 \m<\m \eta \m,$$
so that $y=G(u)=u^2<\eta$, for all $u\in U_\delta\subset[0,\infty)$.

To check Assumption \ref{ass_1}, we use the
result from Remark \ref{rem:Jacob}. In particular, we have that
$$ A(u_0) \m=\m \left. \frac{\partial f(x,u)}{\partial x} \right|_
   {\nm\begin{array}{c} \scriptstyle x=\Xi(u_0)\vspace{-1mm}\\
   \scriptstyle u=u_0 \end{array}} = \begin{bmatrix} -1 & 0 \\ 2u_0 
   & -1 \end{bmatrix}, \vspace{-1mm}$$
for all $u_0\in\rline$. Clearly $\max\Re\sigma(A(u_0))=-1<0$, for all
$u_0\in U_\delta$. Therefore, according to Remark \ref{rem:Jacob}, 
Assumption \ref{ass_1} is satisfied. Clearly $G$ satisfies Assumption
\ref{Ass2}, with $\mu=2(u_{min}-\delta)$. Note that $Y=[u_{min}^2,
u_{max}^2]$. In the following we show that also Assumption 
\ref{ass_3} is met. Solving \rfb{eq:P1}, we get \vspace{-4mm}
$$ \begin{aligned}
   x_1(t) \m=\m e^{-t}x_1(0)+\int_0^te^{\sigma-t}u(\sigma)\dd\sigma,\\
   x_2(t) \m=\m e^{-t}x_2(0)+\int_0^te^{\sigma-t}x_1^2(\sigma)
          \dd\sigma.\\ \end{aligned}$$
Shifting the trajectories of $x_1(t),x_2(t)$ respectively by
$\Xi_1(u_0)=u_0$ and $\Xi_2(u_0)= u_0^2$, we have \vspace{-2mm}
$$ \begin{aligned} x_1(t)-u_0 \m=\m e^{-t}\left[x_1(0)-u_0\right]+
   \int_0^te^{\sigma-t}\left[u(\sigma)-u_0\right] \dd\sigma, \\
   x_2(t)-u_0^2 = e^{-t}\left[x_2(0)-u_0^2\right] + \int_0^t 
   e^{\sigma-t}\left[x_1^2(\sigma)-u_0^2\right] \dd\sigma.
   \end{aligned}$$
Assume that $\abs{u(t)-u_0}\leq m$ for all $t\geq0$, $m>0$. Then
$$\abs{x_1(t)-u_0}\leq e^{-t}\abs{x_1(0)-u_0}+m \FORALL t\geq 0.$$
From the above, it follows trivially that
\begin{equation} \label{eq:lsx_1}
   \limsup_{t\rarrow\infty} \m\abs{x_1(t)-u_0} \m\leq\m m.
\end{equation}
Denote $v(t)\vcentcolon=x_1^2(t)-u_0^2$, then
\begin{multline*}
   \abs{v(t)}=\abs{x^2_1(\sigma)-u_0^2} \m\leq\m \abs{x_1(\sigma)
   -u_0}\cdot\big[\abs{x_1(\sigma)-u_0}\\
   +2(u_{max}+\delta)\big] \leq m\left[m+2(u_{max}+\delta)\right].
\end{multline*}
By a similar reasoning that has lead to \rfb{eq:lsx_1}, we have
$$ \limsup_{t\rarrow\infty} \m\abs{x_2(t)-u_0^2} \m\leq\m m\left[
   m+2(u_{max}+\delta)\right].$$
Combining this with \rfb{eq:lsx_1}, we see that
$$ \limsup_{t\rarrow\infty}\norm{x(t)-\Xi(u_0)}\leq m\left[m+
   2(u_{max}+\delta)+1\right]=\vcentcolon\gamma(m)$$
and this implies Assumption \ref{ass_3}.

\subsubsection{Linearized system stability analysis}

Choose $u_0\in U_\delta$. Linearizing $\PPP$ around the equilibrium
point $(u_0,\Xi(u_0))$, we get the transfer function \vspace{-2mm}
\begin{equation} \label{eq:tf_P1}
   P(s) \m=\m \frac{2u_0}{(s+1)^2}.
\end{equation}
The linearized closed-loop system is given by \rfb{eq:tf_P1} in 
feedback with \rfb{eq:C}. The denominator of its sensitivity is
$$ s^3+2s^2+(1+2u_0k\tau_p)s+2u_0k,$$
where $k>0$ and $\tau_p\geq 0$. Using the well-known Routh test, 
assuming that $\tau_p\in[0,\half)$, the closed-loop is stable iff
$$ k<\frac{1}{u_0(1-2\tau_p)}.$$
Finally, in order to guarantee the closed-loop stability of the
linearized system for all $u_0\in U$, we must have
\begin{equation}\label{eq:kappa_lin}
	k \leq \kappa_{lin}= \frac{1}{(u_{max}+\delta)(1-2\tau_p)}.
\end{equation}

\subsubsection{Nonlinear system stability analysis} 

Substituting \rfb{eq:P1} in \rfb{Trump}, we get the following
closed-loop system equations:
\begin{equation} \label{Bolton}
   \begin{gathered} \dot{x}_1 = -x_1+u_I+\tau_pk(r-{\mathrm B}_\eta(x_2)), \quad 
   \dot{x}_2 = x_1^2-x_2, \\ \dot{u}_I \m=\m \Sscr(u_I,k(r-{\mathrm B}_\eta(x_2))),
   \end{gathered}
\end{equation}
where $r\in Y$ is a constant reference. Using the notation of 
\rfb{eq:P}, \rfb{eq:beta}, and \rfb{eq:x_tilde_u_tilde}, it is clear 
that
$$\begin{gathered} f(x,u_I) \m= \begin{bmatrix} -x_1+u_I \\
   x_1^2-x_2 \end{bmatrix}, \quad g(x) =\m x_2, \\
   \beta(x,u_I,k) \m= \begin{bmatrix} \tau_pk(r-{\mathrm B}_\eta(x_2)) \\ 0 
   \end{bmatrix}, \quad u_r=\m \sqrt{r}, \quad x_r=
   \begin{bmatrix} \sqrt{r} \\ r \end{bmatrix}. \end{gathered}$$
The system \rfb{Bolton} can then be formulated as \rfb{Trump}. We can
apply Theorem \ref{thm:ISS_stability}, yielding that for $k\leq
\kappa^*$ the closed-loop system is globally asymptotically stable and
the tracking error tends to zero. Note that $\kappa^*$ is, in general,
less that $\kappa_{lin}$ from \rfb{eq:kappa_lin}.

\begin{figure} 
   \centering\includegraphics[width=0.5\textwidth]{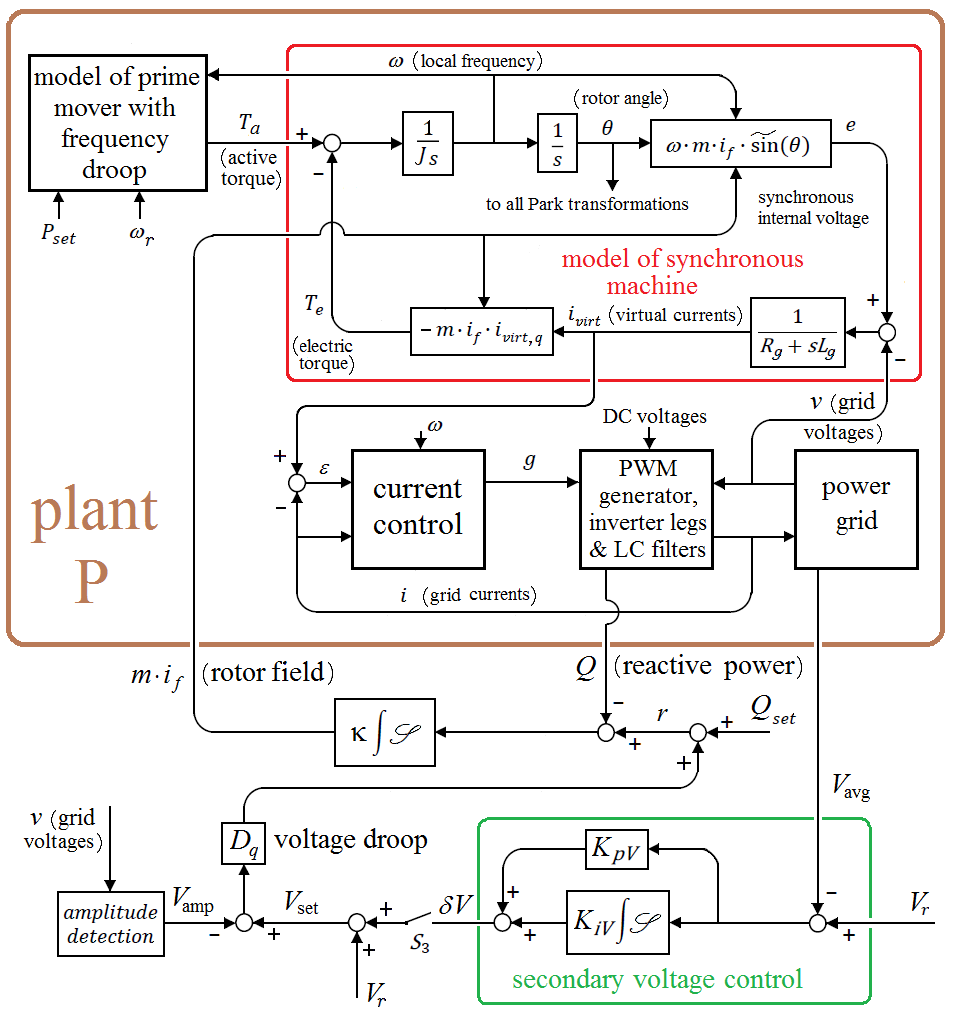}
   \caption{A synchronverter with primary and secondary voltage
   control loops, based on Fig. 2 from \cite{Kustanovich2020b}
   (rearranged). The central secondary voltage control sends the
   signal $\delta V$ to each synchronverter in the microgrid. $V_{\rm
   avg}$ is the (weighted) average of several inverter voltage
   amplitudes.} \label{huge_block_diag}
\end{figure} \vspace{1mm}

\subsection{Voltage control for an inverter} \label{ex:2}

Here we consider the primary and secondary voltage control loops for a
grid-connected synchronverter. For the background on synchronverters
(a special type of 3 phase inverters that emulate synchronous
machines) we refer to \cite{ZhWe:11,Natarajan2017,Kustanovich2020b},
and the references therein. The primary
voltage control loop in a synchronverter is meant to react rapidly (on
the time scale of tens of milliseconds) to changes in the amplitude
$V_{\rm amp}$ of the AC terminal voltages, adjusting the field current
$i_f$ (and hence the field $mi_f$) of the virtual rotor, so that the
reactive power $Q$ flowing from the inverter to the grid tracks the
expression
$$r \m=\m Q_{\rm set} +D_q(V_{\rm set}-V_{\rm amp}) \m.$$
This means that in steady state, $r-Q$ tends to zero. Here $Q_{\rm set}$
is the set value for the reactive power, $D_q$ is the voltage droop
coefficient and $V_{\rm set}$ is the desired amplitude of the terminal
voltage. This primary voltage control loop can be seen in Fig.
\ref{huge_block_diag}, the controller (a saturating integrator with
gain $\kappa$) is just below the large brown block $\PPP$ that
represents the inverter, the power grid and the other blocks of the
synchronverter algorithm. The diagram is based in Figure 2 in
\cite{Kustanovich2020b}, and we refer the reader to the cited
reference for the details of the diagram and the meaning of various
variables that appear there, although we have strived to make the
picture directly understandable to a reader familiar with inverters.

When several inverters are connected into an islanded microgrid, then
it is useful to have a central controller for secondary voltage
control. The purpose of the secondary control is to keep the voltages
in the microgrid within reasonable bounds, by imposing a certain
reference value $V_r$ for a weighted average $V_{\rm ave}$ of the
terminal voltages amplitudes $V_{\rm amp}$ of the individual
inverters. Assuming that the line impedances are low, if in steady
state $V_{\rm avg}=V_r$, then all the individual terminal voltages
will be within reasonable bounds. The secondary voltage controller
proposed in \cite{Kustanovich2020b} is a PI controller with saturating
integrator, as discussed in this paper. Its output is a correction
$\delta V$ that is added to the voltage reference $V_r$ of each
inverter. The secondary controller can be seen on the bottom of
Fig. \ref{huge_block_diag}, and this figure shows clearly that the two
controllers are in a nested structure.

\begin{figure} \centering 
	\includegraphics[width=0.45\textwidth]{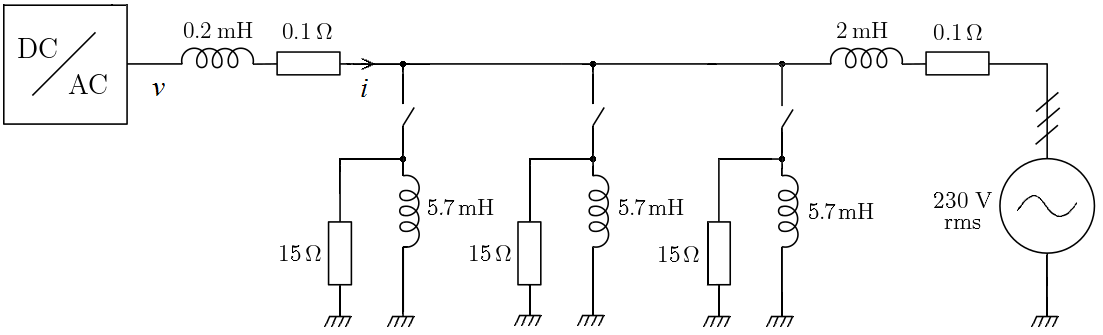}
	\caption{The electric circuit of the system considered
		in Example \ref{ex:2}.} \label{big_load}
\end{figure}

\begin{figure} \centering 
	\includegraphics[width=0.5\textwidth]{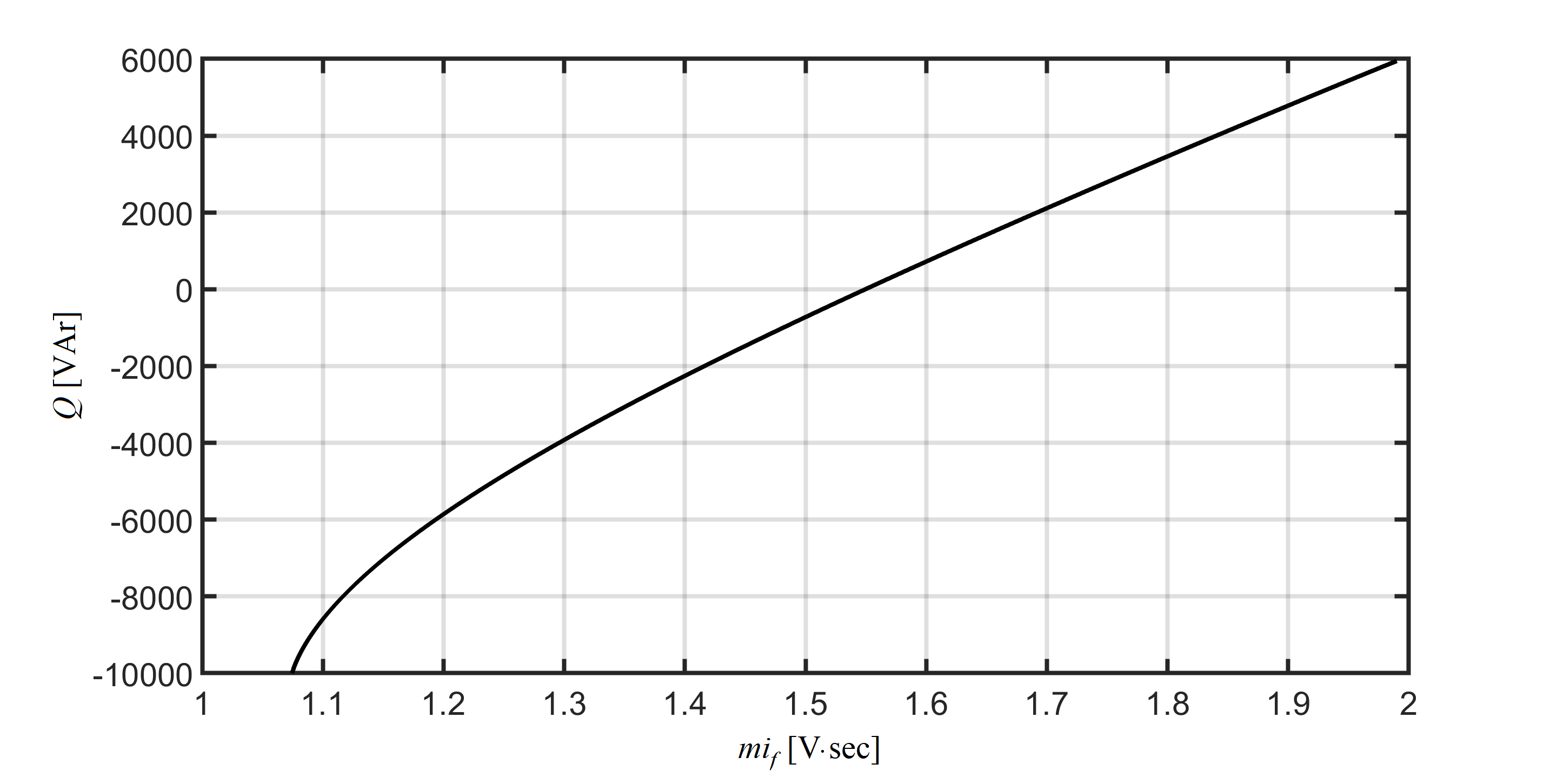}
	\caption{$Q$ as a function of $mi_f$, at steady state,
		for various values of constant $mi_f$, in Example
		\ref{ex:2}.} \label{monotone_plot}
\end{figure}

\begin{figure*}[h!] 
   {\subfigure[Using classical integrators: After the large load is
   connected (in three steps) approximately at $t=10$\m sec, the 
   current increases beyond the allowed range (marked in dashed red).]
   {\includegraphics[width=0.45\textwidth]{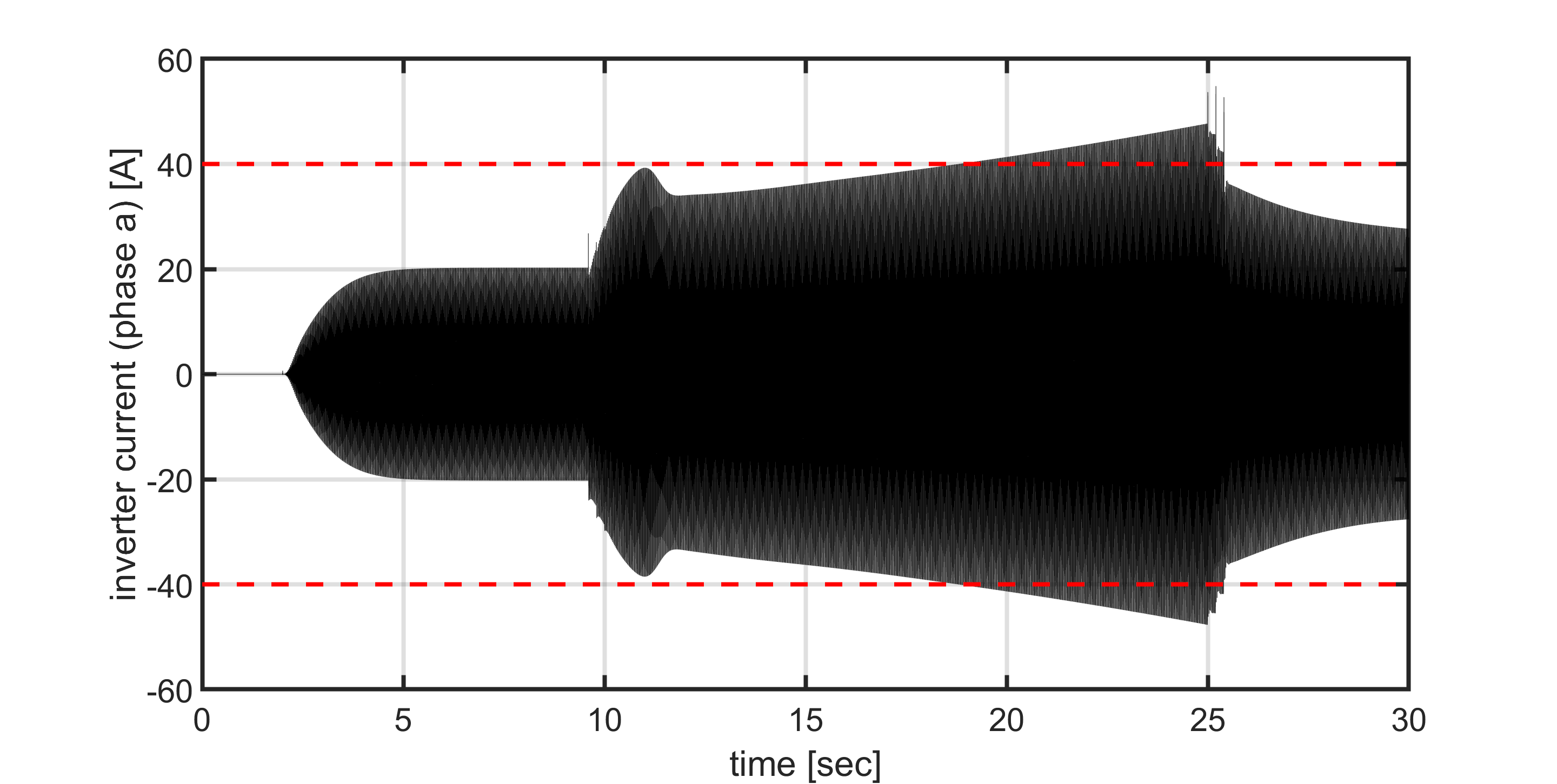}
   \label{plot_no_sat}} \hfil \subfigure[Using saturating integrators:
   After the large load is connected at $t=10$\m sec, the current
   stablizes quickly to a higher amplitude within the allowed range
   (marked in dashed red).]
   {\includegraphics[width=0.45\textwidth]{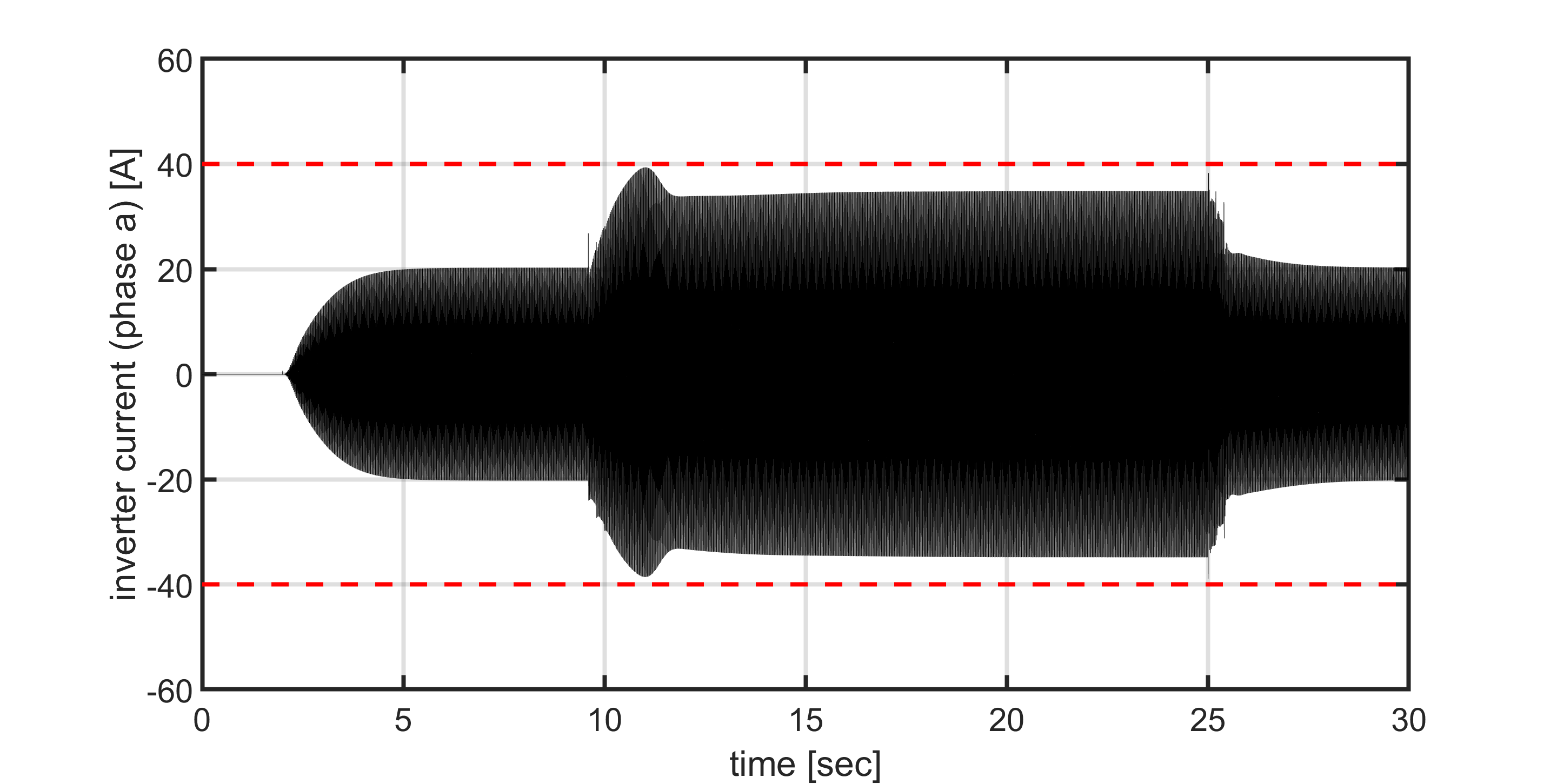}
   \label{plot_with_sat}}} \caption{Connecting the large consumer
   causes a voltage drop from about 10\m sec to about 25\m sec, leading to
   an increase of the inverter current. The plots show the inverter
   current in phase a when the I and PI controllers on the bottom of 
   Fig. \ref{huge_block_diag} are without saturation (Subfig. a) and 
   with saturation (Subfig. b).} \label{fig_sim}
\end{figure*}

For our simulation experiment, we have assumed that a single inverter
is connected to a power grid representing a large number of other
inverters, so that the other inverters together act as an AC voltage
source (also called infinite bus), with frequency 50\m Hz. There is a
line impedance consisting of a 2\m mH inductor and a 0.1\m $\Omega$
resistor from the inverter to the grid, as shown in Fig.
\ref{big_load}. There is a large consumer that connects directly to
the inverter for a period of 15\m sec, starting at about 10\m sec. (The
load is split into 3 equal parallel parts to avoid the large current
peaks that occur with the connection or disconnection of a large
load. These parts act with a delay of 200\m msec between them, so that
it takes 400\m msec to completely connect or to disconnect the load.)
During the first 2\m sec of the simulations, the synchronverter goes
through the start-up procedure and then it connects to the grid. At
about $t=5$\m sec the grid currents reach steady state, with an
amplitude of about 20\m A (the nominal current of this inverter), as can
be seen in Fig. \ref{plot_no_sat} and \ref{plot_with_sat}. In our
simulation experiment, it is assumed that $V_{\rm amp}$ contributes
1/3 of $V_{\rm avg}$, the rest being the infinite bus voltage
amplitude, 325.3\m V.

For lack of space, we do not provide all the details of the analysis
of these control loops. A detailed stability analysis of a closely
related but simpler system appears in \cite{Lorenzetti2020a}, where
the inverter is connected directly to the infinite bus, without all
the other circuit elements in Fig. \ref{big_load}, and there is no 
secondary voltage control ($\delta V=0$ in Fig. 
\ref{huge_block_diag}).

In order for the theory presented here to be applicable to the primary
control loop, we need that the plant $\PPP$ should be described by
sufficiently smooth functions as in \rfb{eq:P}, it should be
exponentially stable around the equilibrium points corresponding to
constant inputs $m i_f$ in a certain range $mi_f\in[u_{min},u_{max}]$
and in steady state, the dependence of the output $Q$ on $mi_f$ (the
function $G$ from Assumption 1) should be monotone increasing. The
smoothness property is easy to verify if we use the average model for
the inverter legs (that are switching at a frequency of several times
10\m kHz). The exponential stability for constant $mi_f$ is a delicate
property that has been proved, under certain sufficient conditions on
the parameters, in \cite{Barabanov2017} and in \cite{Natarajan2018}
(the conditions that they found are not equivalent). Finally, we have
verified the monotone property numerically for a set of typical
parameters for a 10 kW synchronverter, the same parameters as used in
\cite{Kustanovich2020b}. We refrain from listing the parameters here,
but we show the plot of $Q$ (in steady state) as a function of $mi_f$
in Figure \ref{monotone_plot}. A similar plot of monotonicity can be
obtained for the secondary control loop (where the input of the new
plant is $\delta V$ and its output is $V_{\rm avg}$), but we omit the
plot.

During the time when the large load is connected, the terminal voltage 
of the synchronverter drops considerably, causing the signal $V_{\rm 
avg}$ to drop as well, so that $\delta V$ starts growing. If the 
two PI controllers under discussion have normal integrators, then 
$\delta V$ keeps increasing during the time of low voltage and pushes
the signal $mi_f$ to increase as well. This increases the synchronous 
internal voltage $e$, hence the currents, until the currents reach
almost 50\m A towards the end of the low voltage period, see Fig. 
\ref{plot_no_sat}. In reality, this cannot happen, because when the 
currents reach 40\m A (double the nominal current), the protection logic
of the inverter disables the gate signals and opens the circuit
breaker.

By contrast, if we do the same experiment but using two saturating
integrators in the two nested controllers, then the current stabilizes
at around 34\m A during the low voltage period, and quickly returns to
normal steady state after the large consumer is disconnected at around
$t=25$\m sec, see Fig. \ref{plot_with_sat}. Notice that the recovery to
steady state currents is much quicker in Fig. \ref{plot_with_sat} than
in Fig. \ref{plot_no_sat}. 

These different behaviors are due to the secondary voltage control
integrator state, which, in the case of a classical integrator, tends
to windup to excessively large values.

\section{Conclusion} 

We have presented a PI controller with anti-windup (saturating
integrator), for a single-input single-output stable nonlinear plant,
based on singular perturbations theory. Under reasonable assumptions,
we proved that the closed-loop system from Fig. \ref{fig:clos} is able
to track a constant reference signal $r$, while not allowing the
integrator state $u_I$ to leave the set $U$. Moreover, the equilibrium
point of the closed-loop system is stable and has a ``large'' region
of attraction. If the plant is ISS, then the closed-loop system is
globally asymptotically stable. Finally, the validity of our
theoretical results has been demonstrated through two examples: a toy
example and a real application in the control of a certain type of
three phase inverters.

\appendices

\section{Proof of Proposition \ref{prop:uniqueness}} \label{app:A}

\begin{figure} \centering 
   \includegraphics[width=0.4\textwidth]{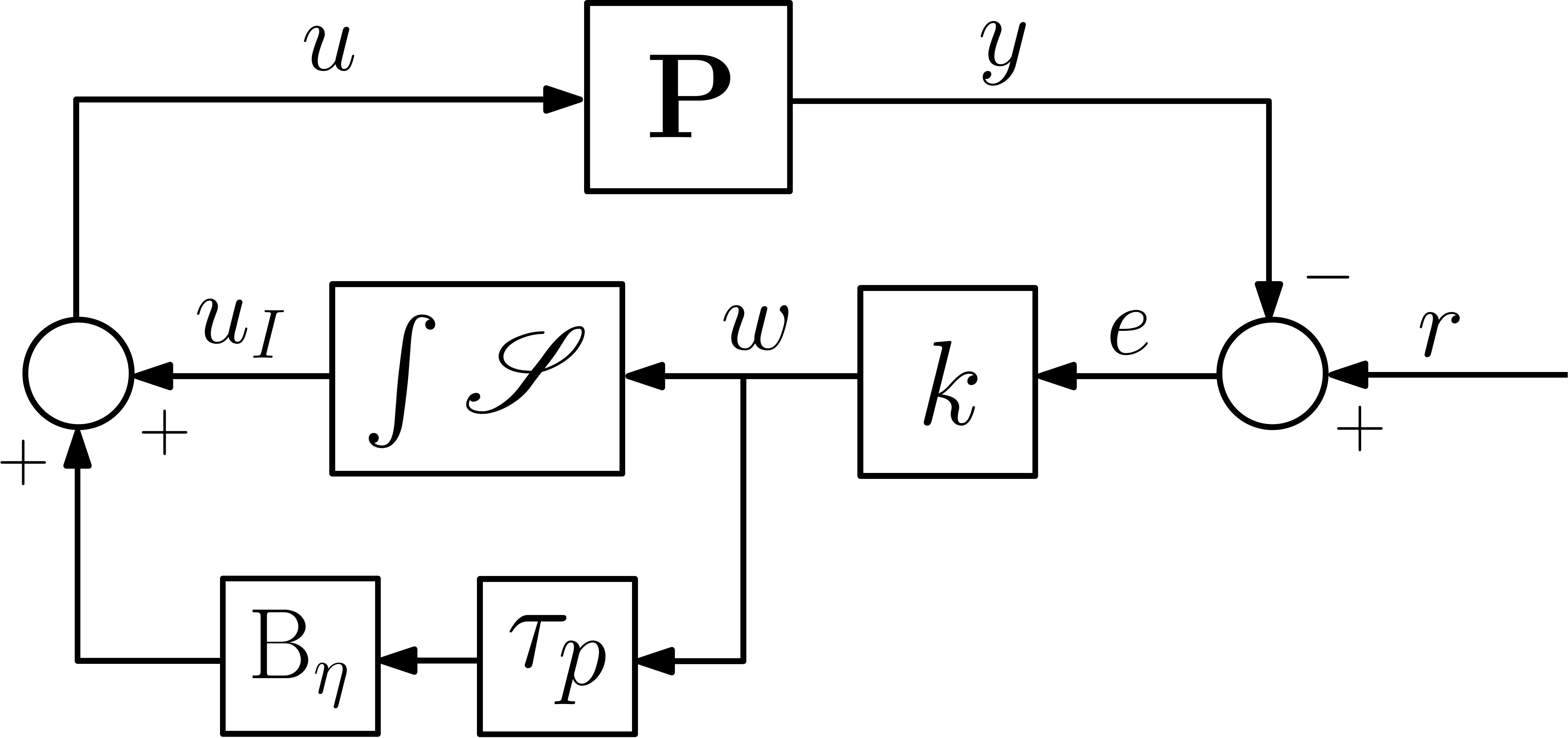}
   \caption{The modified closed-loop system formed by the plant 
   $\PPP$, the saturating integrator $\int\Sscr$, the constant 
   gains $k>0$, $\tau_p\geq 0$ and the saturation function 
   $\mathrm{B}_\eta$, with the reference signal $r$.} \label{fig:clos_eta}
\end{figure}

For the sake of the proof, we slightly modify the closed-loop system
in Fig. \ref{fig:clos} by adding a saturation block, using the
function $\mathrm{B}_\eta $ defined in \rfb{eq:B_eta}, where $\eta>0$
is fixed, to be determined later. The corresponding modified
closed-loop system is shown in Fig. \ref{fig:clos_eta}, so that
$\mathrm{B}_\eta$ is applied to $\tau_p w$, and it is described by
\begin{equation} \label{eq:sys_eta}
   \begin{cases} \begin{gathered}
   \m \ \dot{x} \m=\m f(x,u_I+\mathrm{B}_\eta(\tau_p k(r-g(x)))), \\
   \dot{u}_I \m=\m \Sscr(u_I,k(r-g(x))). \end{gathered} \end{cases}
\end{equation}

\textbf{Notation.} For any $R>0$, $B_R$ is the closed ball of
radius $R$ in $\mathbb{R}^n$. For any $\tau>0$ we denote by $C_\tau$
the set of all the continuous functions on the interval $[0,\tau]$,
with values in $U_\eta\vcentcolon=[u_{min}-\eta,u_{max}+\eta]$. 
This is a complete metric space with the distance
induced by the supremum norm $\norm{\cdot}_\infty$.

{\color{blue}
\begin{lemma} \label{lmm:T_Lips}
Take $x_0\in\rline^n$ and $R>0$. Denote $M=\max\{\|f(x,v)\| \
|\ x\in x_0+B_R, \ v\in U_\eta\}$ and let $\tau\in(0,R/M]$. Then, 
for every $u\in C_\tau$, \rfb{eq:P} has a unique solution $x\in C^1
([0,\tau];\rline^n)$ and $x(t)\in x_0+B_R$ for all $t\in[0,\tau]$.
		
The (nonlinear) operator $T_\tau$ determined by $\PPP$, that maps
any input function $u\in C_\tau$ into an output function $y\in
C[0,\tau]$ (corresponding to $x(0)=x_0$) is Lipschitz continuous,
with Lipschitz constant of the form $\tau L_T$, where $L_T>0$ is
independent of $\tau$.
\end{lemma}}

\begin{IEEEproof} 
It follows from classical ODE theory and the mean value theorem that
for any $u\in C_\tau$, the state trajectory $x(t)$ of $\PPP$ exists
and remains in $x_0+B_R$ for all $t\leq \tau$. Let $L_1,L_2>0$ be such
that for any $z_1,z_2\in x_0+B_R$, and for any $v_1,v_2\in U_\eta$
\begin{equation} \label{eq:Lips_L_1_L_2} 
   \|f(z_2,v_2)-f(z_1,v_1)\| \m\leq\m L_1\|z_2-z_1\|+L_2\|v_2-v_1\|.
\end{equation} 
Such $L_1,L_2$ exist since $f\in C^2$ and $\{x_0+B_R\}\times U_\eta$
is compact. Take two state trajectories of $\PPP$, $x_1$ and $x_2$,
starting from the same initial state $x_0$ and corresponding to inputs
$u_1$ and $u_2$. Then
\begin{displaymath} x_2(t)-x_1(t)=\int_0^t \big[f(x_2(\sigma),
   u_2(\sigma))-f(x_1(\sigma),u_1(\sigma))\big] \dd\sigma,
\end{displaymath} 
for any $t\in[0,\tau]$, whence, using \rfb{eq:Lips_L_1_L_2},
\begin{displaymath} \norm{x_2(t)-x_1(t)} \leq L_1\int_{0}^{t}
   \norm{x_2(\sigma) -x_1(\sigma)} \dd\sigma+\tau L_2
   \norm{u_2-u_1}_\infty.  
\end{displaymath}
It follows from Gronwall's inequality that
\begin{displaymath} \norm{x_2(t)-x_1(t)}\leq\tau L_2\norm{u_1-
   u_2}_\infty\;e^{L_1t} \FORALL t\in[0,\tau],
\end{displaymath}
which implies that
$$ \|x_2-x_1\|_\infty \m\leq\m \tau L_2\|u_1-u_2\|_\infty\;
   e^{L_1\tau}.$$
Finally, since $g\in C^1$, we get
$$\|y_2-y_1\|_\infty \m\leq\m \tau L_T^\tau \|u_1-u_2\|_\infty,$$
where $L_T^\tau=L_2e^{L_1\tau}L_g$ and $L_g$ is the Lipschitz constant
of $g$ on the set $x_0+B_R$. Taking the largest possible value for
$\tau$, i.e., $\tau=\frac{R}{M}$, we get $L_T=L_2e^{L_1\frac{R}{M}}
L_g$. \end{IEEEproof}

\textit{Proof of Proposition \ref{prop:uniqueness}:} We use the
notation from Lemma \ref{lmm:T_Lips}, in particular, $R>0$ is
arbitrary and $x_0$, $u_0$ are fixed. As we have shown in the previous
lemma, the Lipschitz bound of $T_\tau$ is $\tau L_T$, where $L_T>0$ is
independent of $\tau\in(0,R/M]$. We choose $\eta>\delta$ large enough
such that $\abs{\tau_pke(0)}<\eta$, where $e(0)=r-g(x_0)$. Let us
denote by $S_\tau$ the input to output map of the saturating
integrator on the time interval $[0,\tau]$, corresponding to
$u_I(0)=u_0$, and let $\mathrm{B}_\eta$ denote the (Lipschitz
continuous) operator on $C[0,\tau]$ obtained from the function
$\mathrm{B}_\eta$ by applying it pointwise. The estimate
\rfb{eq:lipschitz_u_w} shows that $S_\tau:C[0,\tau]\mapsto C_\tau$ is
Lipschitz continuous with Lipschitz constant $\tau$. If $(x,u_I)$ is a
state trajectory of the modified closed-loop system which is defined
on $[0,\tau]$, and we define $e=r-g(x)$, then we must have (see
Fig. \ref{fig:clos_eta})
$$e \m=\m r - T_\tau (S_\tau + \mathrm{B}_\eta\tau_p) k \m e.$$ 
This can be regarded as a fixed point equation on $C[0,\tau]$. For
$\tau$ sufficiently small so that $\tau L_T\cdot(\tau+\abs{\tau_p})\cdot
\abs{k}<1$, the above equation has a unique solution $e$ according to the
Banach fixed point theorem, see for instance \cite[Sect.~3]{Brooks2009}.
It is easy to see that if $e$ is a solution of the
fixed point equation, if $u=(S_\tau+\mathrm{B}_\eta\tau_p)ke$, $x$ is
the corresponding state trajectory of $\PPP$ starting from $x_0$,
and $u_I=S_\tau ke$, then $(x,u_I)$ is the desired unique state
trajectory of \rfb{eq:sys_eta} on $[0,\tau]$.
	
To extend the above result to the closed-loop system shown in
Fig. \ref{fig:clos}, recall from the beginning of the proof that
$\abs{\tau_pke(0)}<\eta$, i.e., the saturation function $\mathrm{B}
_\eta$ from Fig. \ref{fig:clos_eta} is not saturated at $t=0$. Since
$e(t)=r-g(x(t))$ is a continuous function, there exists $\tau^\prime
\in(0,\tau]$ such that $\abs{\tau_pke(t)}<\eta$ remains true for $t
\in[0,\tau^\prime)$. Then, on the time interval $[0,\tau^\prime)$ the
equations satisfied by $(x,u_I)$ are precisely \rfb{Trump}. Hence,
$(x,u_I)$ is the unique solution of \rfb{Trump} on $[0,\tau^\prime)$.
This $\tau^\prime$ is what we have called $\tau$ in the proposition.
	
Finally, using standard arguments, see for instance Exercise 3.26 in
\cite{Khalil2002}, we get that if an interval of existence $[0,\tau)$
is such that $\tau>0$ is finite and maximal, then 
$\limsup_{t\rarrow\tau}\|(x(t),u_I(t))\|=\infty$,
implying that $\limsup_{t\rarrow\tau}\|x(t)\|=\infty$. \hfill \IEEEQED

\begin{remark}
An alternative way to prove existence (but not uniqueness) of the
solution to the closed-loop system \rfb{Trump} is to use tools from
\textit{differential inclusions} theory. In particular, \rfb{eq:sat}
can be regarded as a \textit{constrained differential inclusion}, see
for example (5.4) in \cite{Goebel2012}, where the state $u$ is
constrained in the closed set $U$ and the map $\Sscr$ is replaced with
the Krasovskii set-valued map \vspace{-2mm}
\begin{displaymath}
   \begin{gathered} \Sscr_K(u_I,w) \m=\m \begin{cases} \{w^+\} &
   \text{if} \quad u_I<u_{min}, \\ [w,w^+] & \text{if} \quad u_I=
   u_{min},\\ \{w\} & \text{if} \quad u_I \in (u_{min},u_{max}), \\ 
   [w^-,w] & \text{if} \quad u_I=u_{max}, \\ \{w^-\} & \text{if} 
   \quad u_I>u_{max}. \end{cases} \end{gathered} 
\end{displaymath} 
The set-valued map describing the closed-loop system is \textit{outer
semicontinuous} (see Definition 5.9 of \cite{Goebel2012}). Therefore,
the closed-loop system \rfb{Trump} satisfies Assumption 6.5 of
\cite{Goebel2012} and, according to Theorem 6.30 of \cite{Goebel2012},
it is \textit{well-posed}.
\end{remark}

\section{Proof of Theorem \ref{thm:11.4}} \label{app:B} 

By assumption, there exist two Lyapunov functions $V$ and $W$, for
the reduced system \rfb{eq:reduced_model} and for the boundary-layer
system \rfb{eq:fast_system_k_0}, respectively, so that the conditions
presented in the first 12 lines of the proof of Theorem 11.4 in
\cite{Khalil2002} hold (our \rfb{eq:V_condit}, \rfb{eq:W_lemma_9_8}).
The translation between our notation and the one in \cite{Khalil2002}
is, again, given by Table \ref{tb:corr}. (Note that, unlike in Theorem
11.4 of \cite{Khalil2002}, assuming that the origin of
\rfb{eq:reduced_model} is exponentially stable is not enough in our
case to apply the converse Theorem 4.14 from \cite{Khalil2002}, since
$\tilde{h}\notin C^1$.) If we define $z$ as in
\rfb{eq:variable_change}, then equations (11.49)-(11.50) in
\cite{Khalil2002} translate to our
\rfb{cl_loop_u_tilde}-\rfb{eq:cl_loop_z}.

We prove some estimates. Since the system
\rfb{cl_loop_u_tilde}-\rfb{eq:cl_loop_z} is time-invariant and
$\tilde{h}$ is independent of $k$, the growth conditions (valid in a
neighbourhood of the origin) shown in \cite{Khalil2002} (after
(11.50)) reduce in our case to
\begin{subequations}\label{eq:gc}
   \begin{gather}
   \big\|\tilde{\beta}(\tilde{u}_I,z+\tilde{\Xi}(\tilde{u}_I),k)
   \big\|\leq k\alpha_1(|\tilde{u}_I|+\norm{z}), \label{eq:gc1} \\
   \big|\tilde{h}(\tilde{u}_I,z+\tilde{\Xi}(\tilde{u}_I))-
   \tilde{h}(\tilde{u}_I,\tilde{\Xi}(\tilde{u}_I))\big|\leq
   \alpha_2\norm{z}, \label{eq:gc2} \\
   \big|\tilde{h}(\tilde{u}_I,\tilde{\Xi}(\tilde{u}_I))
   \big|\leq \alpha_3|\tilde{u}_I|, \quad
   \bigg\|\frac{\dd\tilde{\Xi}}{\dd\tilde{u}_I}\bigg\|
   \leq \alpha_4, \label{eq:gc34} \end{gather}
\end{subequations}
where $\alpha_1,\dots \alpha_4$ are positive constants. These
estimates follow respectively from:
\begin{itemize}
   \item $\tilde{\beta}\in {\rm Lip}$ and $\tilde{\beta}(\tilde{u}_I,\tilde
      {x},0)=0$ for all $(\tilde{u}_I,\tilde{x})\in\tilde{X}_\delta$,
   \item $\tilde{h}$ being uniformly Lipschitz in the second
      argument, \item $\tilde{h}$ being uniformly Lipschitz in the
      second argument, $\tilde{\Xi}\in C^1$, $\tilde{h}(0,0)=0$, and
      $\tilde{\Xi}(0)=0$.
   \item $\tilde{\Xi}\in C^1$.
\end{itemize}

Consider the function $\nu$ from \rfb{eq:nu_Lyap} as a Lyapunov 
function candidate for the system 
\rfb{cl_loop_u_tilde}-\rfb{eq:cl_loop_z}. Computing its derivative
along the state trajectories of 
\rfb{cl_loop_u_tilde}-\rfb{eq:cl_loop_z}, we get
\begin{equation} \label{eq:Lyap_cl}
   \frac{\dd\nu}{\dd s} = \frac{\dd V}{\dd \tilde{u}_I}\frac{\dd 
   \tilde{u}_I}{\dd s} + \frac{\partial W}{\partial \tilde{u}_I}
   \frac{\dd \tilde{u}_I}{\dd s}+\frac{\partial W}{\partial z}
   \frac{\dd z}{\dd s}.
\end{equation}
Consider the first term on the right-hand side of \rfb{eq:Lyap_cl}. 
Substituting \rfb{cl_loop_u_tilde} we obtain \vspace{-2mm}
\begin{multline*}
   \frac{\dd V}{\dd \tilde{u}_I}\tilde{h}(\tilde{u}_I,z+\tilde{\Xi}
   (\tilde{u}_I)) \leq\frac{\dd V}{\dd \tilde{u}_I}\tilde{h}
   (\tilde{u}_I,\tilde{\Xi}(\tilde{u}_I)) \\ + \bigg|\frac{\dd V}{\dd 
   \tilde{u}_I}\bigg|\big|\tilde{h}(\tilde{u}_I,z+\tilde{\Xi}
   (\tilde{u}_I))-\tilde{h}(\tilde{u}_I,\tilde{\Xi}(\tilde{u}_I))
   \big|.
\end{multline*}
Using \rfb{eq:V_condit} and \rfb{eq:gc2}, we have
\begin{equation} \label{eq:b1}
   \frac{\dd V}{\dd \tilde{u}_I}\frac{\dd \tilde{u}_I}{\dd s}\leq
   -c_3\abs{\tilde{u}_I}^2+\alpha_2c_4\abs{\tilde{u}_I}\norm{z}.
\end{equation}
Consider the second term on the right-hand side of \rfb{eq:Lyap_cl}.
Substituting \rfb{cl_loop_u_tilde} we obtain
\begin{multline*}
   \frac{\partial W}{\partial \tilde{u}_I}\tilde{h}(\tilde{u}_I,z+
   \tilde{\Xi}(\tilde{u}_I))\leq\bigg\|\frac{\partial W}
   {\partial\tilde{u}_I}\bigg\|\bigg[\big|\tilde{h}(\tilde{u}_I,
   \tilde{\Xi}(\tilde{u}_I))\big|\\+\big|\tilde{h}(\tilde{u}_I,z +
   \tilde{\Xi}(\tilde{u}_I))-\tilde{h}(\tilde{u}_I,\tilde{\Xi}
   (\tilde{u}_I))\big|\bigg].
\end{multline*}
Using \rfb{eq:W_lemma_9_8}, \rfb{eq:gc2} and \rfb{eq:gc34}, we get
\begin{equation} \label{eq:b2}
   \frac{\partial W}{\partial\tilde{u}_I}\frac{\dd\tilde{u}_I}
   {\dd s}\leq b_5\norm{z}^2\big(\alpha_2\norm{z}+\alpha_3
   |\tilde{u}_I|\big).
\end{equation}
Finally, consider the third term on the right-hand side of
\rfb{eq:Lyap_cl}. Substituting \rfb{eq:cl_loop_z} we obtain
\begin{multline*}
   \frac{\partial W}{\partial z}\bigg[\frac{1}{k}\tilde{f}
   (\tilde{u}_I,z+\tilde{\Xi}(\tilde{u}_I)) + \frac{1}{k}
   \tilde{\beta}(\tilde{u}_I,z+\tilde{\Xi}(\tilde{u}_I),k) \\
   -\frac{\dd\tilde{\Xi}}{\dd\tilde{u}_I}\tilde{h}(\tilde{u}_I,z+
   \tilde{\Xi}(\tilde{u}_I))\bigg]
\end{multline*}
\begin{multline*}
   \leq\frac{1}{k}
   \frac{\partial W}{\partial z} \tilde{f}(\tilde{u}_I,z+
   \tilde{\Xi}(\tilde{u}_I))+\bigg\|\frac{\partial W}{\partial z}
   \bigg\| \bigg[\frac{1}{k}\big\|\tilde{\beta}(\tilde{u}_I,z +
   \tilde{\Xi}(\tilde{u}_I),k)\|\\ + \bigg\|\frac{\dd\tilde{\Xi}}
   {\dd\tilde{u}_I}\bigg\|\big\|\tilde{h}(\tilde{u}_I,z+\tilde{\Xi}
   (\tilde{u}_I))\big\|\bigg].
\end{multline*}
Proceeding as before and using \rfb{eq:W_lemma_9_8}, \rfb{eq:gc},
we get
\begin{multline} \label{eq:b3}
   \frac{\partial W}{\partial z}\frac{\dd z}{\dd s}\leq-\frac{b_3}
   {k}\norm{z}^2+\alpha_1b_4\norm{z}(\abs{\tilde{u}_I}+\norm{z})\\
   + \alpha_4b_4\norm{z}(\alpha_2\norm{z}+\alpha_3\abs{\tilde{u}_I}).
\end{multline}
Putting together \rfb{eq:b1}, \rfb{eq:b2} and \rfb{eq:b3}, we have from \rfb{eq:Lyap_cl}
\begin{multline*}
   \frac{\dd\nu}{\dd s}\leq-c_3\abs{\tilde{u}_I}^2-\frac{b_3}{k}
   \norm{z}^2+d_1\norm{z}^2 \\ +d_2\abs{\tilde{u}_I}\norm{z}+d_3
   \abs{\tilde{u}_I}\norm{z}^2+d_4\norm{z}^3,
\end{multline*}
with positive $c_3$, $b_3$ and non negative $d_1,\dots d_4$. For all
$z\leq\e_0$, this reduces to
\begin{displaymath}
   \frac{\dd\nu}{\dd s}\leq-c_3\abs{\tilde{u}_I}^2-\frac{b_3}{k}
   \norm{z}^2+d_5\norm{z}^2 + 2d_6\abs{\tilde{u}_I}\norm{z},
\end{displaymath}
for non negative $d_5$, $d_6$. The above inequality can be rewritten
in matrix form as
\begin{displaymath}
   \frac{\dd\nu}{\dd s}\leq-\begin{bmatrix} \abs{\tilde{u}_I} \\
   \norm{z} \\ \end{bmatrix}^T \begin{bmatrix} c_3 & -d_6 \\
   -d_6 & \frac{b_3}{k}-d_5 \end{bmatrix} \begin{bmatrix}
   \abs{\tilde{u}_I} \\ \norm{z} \\ \end{bmatrix}.
\end{displaymath}
(Note that, differently from \cite{Khalil2002}, $k$ does not appear in
the upper-left term, since our reduced model \rfb{eq:reduced_model} is
independent of $k$.) Therefore, there exists a $\kappa>0$ such that
for all $0<k<\kappa$, we have $\frac{\dd\nu}{\dd s}\leq-2\gamma\nu$,
for some $\gamma>0$. It follows that
\begin{displaymath}
\nu(\tilde{u}_I(s),z(s))\leq e^{-\gamma s}\nu(\tilde{u}_I(0),z(0)),
\end{displaymath}
and from \rfb{eq:V_condit}, \rfb{eq:W_lemma_9_8}
\begin{equation} \label{eq:exp_bound}
   \begin{bmatrix} \abs{\tilde{u}_I(s)} \\ \norm{z(s)} \\  
   \end{bmatrix}\leq K_1 e^{-\gamma s} \begin{bmatrix}
   \abs{\tilde{u}_I(0)} \\ \norm{z(0)} \\ \end{bmatrix},
\end{equation}
for some $K_1>0$. Recall the change of variables
\rfb{eq:variable_change}, since $\norm{\tilde{\Xi}(\tilde{u}_I)}\leq
L_{\tilde{\Xi}}\abs{\tilde{u}_I}$ ($\tilde{\Xi}\in C^2$ and $\tilde{U}
_\delta$ compact), then from \rfb{eq:exp_bound} it follows that
\begin{displaymath}
   \begin{bmatrix} \abs{\tilde{u}_I(s)} \\ \norm{\tilde{x}(s)} \\  
   \end{bmatrix}\leq K_2 e^{-\gamma s} \begin{bmatrix}
   \abs{\tilde{u}_I(0)} \\ \norm{\tilde{x}(0)} \\ \end{bmatrix},
\end{displaymath}
for some $K_2>0$, completing the proof.
\hfill \IEEEQED

\section{Proof of Proposition \ref{prop:ISS_sol}} \label{app:C}

As assumed in the proposition, $g$ is bounded or $\tau_p=0$. Choose
an initial state $(x_0,u_0)\in X_\delta$ for the closed-loop system
\rfb{Trump}. From Proposition \ref{prop:uniqueness} it follows that
there exists a maximal $\tau\in(0,\infty]$ such that \rfb{Trump} has a
unique state trajectory $(x,u_I)$ defined on $[0,\tau)$. If $\tau$ is
finite, then $\limsup_{t\rarrow\tau}\|x(t)\|=\infty$. This contradicts
Proposition \ref{prop:ISS}. Indeed, recalling the definition of
$\int\mathscr{S}$ from \rfb{eq:sat}, $u_I(t)\in U_\delta$ for all
$t\geq0$. Thus, according to Proposition \ref{prop:ISS}, $\norm{x(t)}$
is bounded for all $t\geq 0$. Therefore, for any $(x_0,u_0)\in
X_\delta$, \rfb{Trump} has a unique global solution for all $k>0$ and
$r\in Y$.

We now prove continuous dependence of the solution on the initial
state.  Take $(x_0,u_0)$ and $(\bar{x}_0,\bar{u}_0)$ to be initial
states in $X_\delta$.  Let $(x,u_I)$ be the state trajectory of
\rfb{Trump} for the initial state $(x_0,u_0)$, and, similarly, let
$(\bar{x},\bar{u}_I)$ be the state trajectory corresponding to
$(\bar{x}_0,\bar{u}_0)$. Denote by $B_R$ the closed ball of radius $R$
in $\mathbb{R}^n$. Let $M>0$ be such that $\abs{u(t)-v}\leq M$ for 
all $t\geq0$, for any trajectory of \rfb{Trump} starting in $X_\delta$,
and for any $v\in U_\delta$. Such an $M$ exists because $g$ is bounded
or $\tau_p=0$. Then from \rfb{eq:ISS} we have, for any $v\in U_\delta$,
\begin{equation*}
   \begin{aligned} \norm{x(t)-\Xi(v)}\leq\alpha_v(\norm{x_0-\Xi(v)},t)
   +\gamma(M),\\ \norm{\bar{x}(t)-\Xi(v)}\leq\alpha_v
   (\norm{\bar{x}_0-\Xi(v)},t)+\gamma(M),\end{aligned}
\end{equation*}
which implies that there exists $\rho>0$ such that $x(t),\bar{x}(t)
\in B_\rho$ for all $t\geq0$. Let $L_1,L_2>0$ be such that
$$ \norm{f(z_2,v_2)-f(z_1,v_1)} \m\leq\m L_1\norm{z_2-z_1}+L_2
   \abs{v_2-v_1},$$
for any $z_1,z_2\in B_\rho$ and for any $v_1,v_2\in U_\delta$. Such 
$L_1,L_2$ exist since $f\in C^2$ and $B_\rho\times U_\delta$ is
compact. It follows that
\begin{multline} \label{eq:bush1}
   \norm{x(t)-\bar{x}(t)} \leq \norm{x_0-\bar{x}_0} + L_1\int_0^t
   \norm{x(\sigma)-\bar{x}(\sigma)} \dd\sigma\\ +L_2\int_0^t 
   \abs{u_I(\sigma)-\bar{u}_I(\sigma)} \dd\sigma.
\end{multline}
Besides, from \rfb{eq:blair} (extended to continuous functions) we 
have
$$ \abs{u_I(t)-\bar{u}_I(t)} \m\leq\m \abs{u_I(0)-\bar{u}_I(0)} +
   \int_0^t\abs{w(\sigma)-\bar{w}(\sigma)} \dd\sigma,$$
where $w(\sigma)=k(r-g(x(\sigma)))$, $\bar{w}(\sigma)=k(r-g(\bar{x}
(\sigma)))$. Let $L_3$ be the Lipschitz constant of $g$ on $B_\rho$,
then
\begin{equation} \label{eq:bush3}
   \abs{w(\sigma)-\bar{w}(\sigma)} \m\leq\m kL_3 \norm{x(\sigma)-
   \bar{x}(\sigma)},
\end{equation}
for all $\sigma\geq 0$ and $x,\bar{x}\in B_\rho$. Adding
\rfb{eq:bush1} to \rfb{eq:bush3}, we get
\begin{multline*}
   \norm{x(t)-\bar{x}(t)}+\abs{u_I(t)-\bar{u}_I(t)} \leq 
   \norm{x_0-\bar{x}_0}+\abs{u_0-\bar{u}_0}\\ +(L_1+kL_3) \int_0^t
   \norm{x(\sigma)-\bar{x}(\sigma)} \dd\sigma + L_2\int_0^t 
   \abs{u_I(\sigma)-\bar{u}_I(\sigma)} \dd\sigma.
\end{multline*}
Denote $L_4=\max\{L_1+kL_3,L_2\}$, then defining $q(t)\vcentcolon=
\norm{x(t)-\bar{x}(t)}+\abs{u_I(t)-\bar{u}_I(t)}$, the last estimate
becomes
$$q(t) \m\leq\m q(0)+L_4\int_0^tq(\sigma) \dd\sigma.$$
Finally, using the Gronwall inequality, we get
$$q(t) \m\leq\m q(0)e^{L_4t}.$$
This shows that at any $t\geq0$, the state $(x(t),u_I(t))$ from 
\rfb{Trump} depends Lipschitz continuously on the initial state 
$(x_0,u_0)$. \hfill \IEEEQED

\section{Proof of Theorem \ref{thm:ISS_stability}} \label{app:D}

{\color{blue}
\begin{lemma} \label{lmm:X_T}
We work under the assumptions of Proposition \ref{prop:ISS_sol}. Let
$T>0$ and let $X_T\subset X_\delta$ be a compact set such that if
$(x_0,u_0)\in X_T$, then the state trajectory $x$ of $\PPP$ starting
from $x(0)=x_0$, with constant input $u_0$, satisfies
$$\norm{x(T)-\Xi(u_0)} \m\leq\m \e_0 \m.$$
Then there exists a $\kappa_T>0$ such that for any $k\in(0,\kappa_T)$,
for any $r\in Y$, if the initial state $(x_0,u_0)$ of \rfb{Trump} is
in $X_T$, then its state trajectory $(x,u_I)$ satisfies
\begin{equation} \label{eq:conv} 
   x(t)\rarrow\Xi(u_r),\qquad u(t)\rarrow u_r,\qquad y(t)\rarrow r,
\end{equation} 
and this convergence is at an exponential rate.
\end{lemma}}

\begin{IEEEproof} Take $r\in Y$. Recall the changes of variables 
\rfb{eq:x_tilde_u_tilde}, \rfb{eq:variable_change} and the state space
$\tilde{X}_\delta=\tilde{U}_\delta\times\rline^n$ of the closed-loop
system \rfb{cl_loop_u_tilde}-\rfb{eq:cl_loop_z}. Define the set
$\tilde{X}_T\subset\tilde{X}_\delta$ as the image of the set $X_T$
through the change of variables just mentioned. Then for any
$(\tilde{u}_0, z_0)\in\tilde{X}_T$, the state trajectory $z$ of the
boundary-layer system \rfb{eq:fast_system_k_0} starting from
$z(0)=z_0$, with fixed input $\tilde{u}_0$, satisfies
$\norm{z(T)}\leq\e_0$. Clearly $\tilde{X}_T$ is compact.
	
The convergence of $z$ to $0$ can be regarded as an exponential one.
Indeed, for $\norm{z_0}\leq\e_0$ this is clear from \rfb{eq:exp_stab}.
For $\norm{z_0}>\e_0$ we invoke a compactness argument: there exist
$\tilde{m}>0$ such that $\norm{z(t)}\leq\tilde{m}e^{-\l t}\norm{z_0}$
for all $(z_0,\tilde{u}_0)\in\tilde{X}_T$ with $\norm{z_0}\geq\e_0$
and for all $t\in[0,T]$. Then for $t>T$ we have, according to
\rfb{eq:exp_stab}, $\norm{z(t)}\leq \tilde{m}e^{-\l(t-T)}\norm{z(T)}
\leq \tilde{m}me^{-\l t}\norm{z_0}$.
	
Applying Lemma 9.8 of \cite{Khalil2002}, we get a Lyapunov function
$W(\tilde{u}_I,z)$ for the boundary-layer system
\rfb{eq:fast_system_k_0} such that \rfb{eq:W_lemma_9_8} holds for all
$(\tilde{u}_I,z)\in\tilde{X}_T$. The constants $b_i$, $i=1,\dots 5$
from \rfb{eq:W_lemma_9_8} in this case are different, since
$\tilde{m}_T$ and $\tilde{\l}_T$ are, in general, different from $m$
and $\l$ of Assumption \ref{ass_1}. From this point on, we apply
Theorem \ref{thm:11.4} and everything proceeds as in the proof of
Theorem \ref{thm:stab}, substituting the inequality \rfb{eq:W_bound}
with
\begin{equation}\label{eq:biden}
   W(\tilde{u}_I,z)\leq b_2\max\{\norm{z}^2 \ | \ z\in\Pi
   \tilde{X}_T\} \ \ \ \forall \m (\tilde{u}_I,z)\in\tilde{X}_T,
\end{equation}
where $\Pi$ denotes the projection onto the second component in the
product $\tilde{U}_\delta\times\rline^n$, and modifying the definition
of $L$ in \rfb{eq:L} by substituting $b_2\e_0^2$ with the right-hand
side of \rfb{eq:biden}.
\end{IEEEproof}

\textit{Proof of Theorem \ref{thm:ISS_stability}:} For any $(x_0,u_0)
\in X$, we denote by $\phi(t,x_0,u_0)$ the solution of \rfb{eq:P}
corresponding to the initial state $x_0$ and the constant input $u_0$.
Then Assumption \ref{ass_1} implies that for $\norm{z_0}\leq\e_0$ we
have
$$ \norm{\phi(t,x_0,u_0)-\Xi(u_0)} \leq me^{-\l t} 
   \norm{x_0-\Xi(u_0)},$$
for all $t\geq0$. For each $T>0$, we define the relatively open set
$\mathring{X}_T\subset X$ as follows:
$$ \mathring{X}_T \nm\m= \{(x_0,u_0)\nm\mm\in\nm\m X \m|\ \norm{
   \phi(T,x_0,u_0)-\Xi(u_0)} \nm<\nm\mm\e_0,\mm\norm{x_0} \nm\mm<
   \nm\mm T\}.$$
Clearly Assumption \ref{ass_3} implies that every $(x_0,u_0)\in
X$ is contained in some of the sets $\mathring{X}_T$. In other
words, the union of all the sets $\mathring{X}_T$ is $X$.

Let us choose an (arbitrary) upper bound $k_{max}>0$ for the gain $k$
in Fig. \ref{fig:clos}. It follows from the assumptions that there
exists an $M>0$ such that for any state trajectory $(x,u_I)$ of
\rfb{Trump} (with any $r\in Y$, any gain $k\in[0,k_{max}]$ and
starting from any initial state in $X$), the function $u=u_I+\tau_p k
(r-g(x))$ satisfies $|u(t)|\leq M$. Denote $R=\gamma(M)+1+\max_{u_0
\in U}\norm{\Xi(u_0)}$ and, as usual, $B_R$ is the closed ball of
radius $R$ in $\rline^n$. We claim that for any $r\in Y$, any $k\in
[0,k_{max}]$ and for any initial state $(x_0,u_0)\in X$, the
closed-loop system \rfb{Trump} has a solution for all $t\geq 0$ and
there exists $T^*>0$ such that $x(t)\in B_R$, for all $t\geq T^*$.
Indeed, if we choose $T^*>0$ such that $\alpha_{u_0}(\norm{x_0-
\Xi(u_0)},T^*)\leq 1$, then from \rfb{eq:ISS},
$$ \norm{x(t)-\Xi(u_0)} \m\leq\m 1+\gamma(M) \FORALL \m t\geq T^*,$$
which implies our previous claim. Since $B_R$ is compact and the sets
\{$\mathring{X}_T\ |\ T>0$\} are an open covering of $B_R$, we can
extract a finite cover. Since $\mathring{X}_T$ is increasing with $T$,
there exists $T>0$ such that $B_R\subset\mathring{X}_T$. Let us denote
by $X_T$ the closure of $\mathring{X}_T$. Then we can apply Lemma
\ref{lmm:X_T} for this $X_T$, to conclude that there exists $\kappa^*
\in(0,k_{max}]$ such that for any $k\in(0,\kappa^*]$, every state
trajectory of \rfb{Trump} that, at some point, is in $X_T$ satisfies
\rfb{eq:conv} and the convergence is at an exponential rate. But we
have seen earlier that every state trajectory of \rfb{Trump} (starting
from any initial state) reaches $B_R$ in some time. Hence (since
$B_R\subset X_T$), we conclude that for any $k\in(0,\kappa^*]$, any
state trajectory of \rfb{Trump} satisfies $\rfb{eq:conv}$ and the
convergence is at an exponential rate. \hfill\IEEEQED

\section*{Acknowledgments}

We thank Shivprasad Shivratri for setting up the simulation experiment
described in Subsect. \ref{ex:2}. Shivprasad is an MSc student
supported by the grant no. 217-11-037 from the Ministry of
Infrastructure and Energy. We also thank Nathanael Skrepek for
useful discussions.


\begin{IEEEbiographynophoto}{Pietro Lorenzetti} is an Early Stage
Researcher within the Marie Curie ITN project ``ConFlex'', who focuses
his research on nonlinear control. He is working under
the supervision of George Weiss, and his co-supervisor is Enrique Zuazua.
Pietro has completed the bachelor degree in ``Computer engineering and
automation'' at Universita Politecnica delle Marche, in Ancona. In
2015 he graduated with honours and he moved to Torino, where he 
enrolled the master degree in ``Mechatronic Engineering'' at 
Politecnico di Torino. In the same year, he also joined the 
double-degree program ``Alta Scuola Politecnica'', a highly selective
joined program between Politecnico di Torino and Politecnico di
Milano. In 2017 he graduated in both Politecnico di Milano and
Politecnico di Torino, with honours. His research interests include
nonlinear systems, nonlinear control, and power system stability.
\end{IEEEbiographynophoto}

\begin{IEEEbiographynophoto}{George Weiss}
received the MEng degree in control engineering from the Polytechnic
Institute of Bucharest, Romania, in 1981, and the Ph.D. degree in
applied mathematics from Weizmann Institute, Rehovot, Israel, in
1989. He was with Brown University, Providence, RI, Virginia Tech,
Blacksburg, VA, Ben-Gurion University, Beer Sheva, Israel, the
University of Exeter, U.K., and Imperial College London, U.K. His
current research interests include distributed parameter systems,
operator semigroups, passive and conservative systems (linear and
nonlinear), power electronics, repetitive control, sampled data
systems, and the grid integration of distributed energy sources.
\end{IEEEbiographynophoto}
\vfill

\end{document}